\documentclass[10pt]{article}
\usepackage[top=2.5cm, bottom=2cm, inner=2.5cm, outer=2.5cm]{geometry}

\usepackage{xcolor}
\usepackage{booktabs, multirow, makecell}
\usepackage{subfigure, bm, comment}
\usepackage{amsmath,amsfonts,amssymb}
\usepackage{fixmath}
\usepackage{algorithm, algpseudocode}
\usepackage[mathcal]{eucal}
\usepackage{pstool}
\usepackage{hyperref}
\usepackage{pgfplots}
\def\x{\mathbold{x}}
\def\A{\mathbf{A}}
\def\B{\mathbf{B}}

\def\OJ{\mathsf{J}}

\def\Q{\mathbf{Q}}
\def\E{\mathbf{E}}

\def\v{\mathbf{v}}

\def\q{\mathbold{q}}
\def\r{\mathbold{r}}

\def\trace{\mathrm{tr}}

\newcommand{\uu}{\mathbold{u}}

\newcommand{\z}{\mathbold{z}}
\newcommand{\w}{\mathbold{w}}
\def\y{\mathbold{y}}
\def\e{\mathbold{e}}
\def\W{\mathbf{W}}
\def\H{\mathbf{H}}
\def\F{\mathbf{F}}

\def\R{\mathbf{R}}
\def\P{\mathbf{P}}
\def\K{\mathbf{K}}

\newcommand{\reals}{\mathbb R}
\newcommand{\naturals}{\mathbb N}
\newcommand{\prox}{\mathsf{prox}}

\def\X{\mathbf{X}}\def\Y{\mathbf{Y}}
\def\Xcb{{\color{blue}\mathbf{X}}}
\def\Ycb{{\color{blue}\mathbf{Y}}}
\def\Lambdacb{{\color{blue}\mathbf{\Lambda}}}
\def\I{\mathbf{I}}
\def\OI{\mathsf{I}}

\def\Phib{\mathbf{\Phi}}
\def\Psib{\mathbf{\Psi}}
\def\Sigmab{\mathbf{\Sigma}}

\def\transp{\mathsf{T}}
\def\rhop{\rho_{\mathrm{p}}}
\def\rhoc{\rho_{\mathrm{c}}}

\def\gammacb{{\color{blue} \gamma}}

\def\Wcb{{\color{blue} \W}}
\def\lambdacb{{\color{blue} \lambda}}
\def\sigmac{\sigma_{\mathrm{c}}}
\def\zetap{\zeta_{P, \rhop}}
\def\zetac{\zeta_{C, \rhoc}}
\def\xip{\xi_{P, \rhop}}

\def\qed{\hfill $\blacksquare$ }
\DeclareMathOperator*{\minimize}{\textrm{minimize}}

\usepackage{theorem}
\newtheorem{remark}{Remark}
\newtheorem{lemma}{Lemma}
\newtheorem{proposition}{Proposition}
\newtheorem{assumption}{Assumption}
\newtheorem{definition}{Definition}
\newtheorem{theorem}{Theorem}

\newtheorem{example}{\bf Example}
\newenvironment{proof}{\paragraph{\emph{Proof:}}}{\hfill$\diamondsuit$}

\newcommand{\rev}[1]{{\color{black} #1}}

\usepackage{hyperref}

\makeatletter
\def\@fnsymbol#1{\ensuremath{\ifcase#1\or a\or b\or
   \mathsection\or \mathparagraph\or \|\or **\or \dagger\dagger
   \or \ddagger\ddagger \else\@ctrerr\fi}}
    \makeatother

\title{\sf Nonlinear Optimization Filters for Stochastic Time-Varying Convex Optimization}

\author{Andrea Simonetto \thanks{
UMA, ENSTA Paris, Institut Polytechnique de Paris, 91120 Palaiseau, France; \texttt{andrea.simonetto@ensta-paris.fr}. }\qquad
Paolo Massioni \thanks{Univ Lyon, INSA Lyon, Universit\'{e} Claude Bernard Lyon 1, Ecole Centrale de Lyon, CNRS, Amp\`{e}re, UMR5505, 69621 Villeurbanne, France; \texttt{paolo.massioni@insa-lyon.fr}}
}

\begin{document}

\maketitle

\begin{abstract}{We look at a stochastic time-varying optimization problem and we formulate online algorithms to find and track its optimizers in expectation. The algorithms are derived from the intuition that standard prediction and correction steps can be seen as a nonlinear dynamical system and a measurement equation, respectively, yielding the notion of nonlinear filter design. The optimization algorithms are then based on an extended Kalman filter in the unconstrained case, and on a bilinear matrix inequality condition in the constrained case. Some special cases and variations are discussed, notably the case of parametric filters, yielding certificates based on LPV analysis and, if one wishes, matrix sum-of-squares relaxations. Supporting numerical results are presented from real data sets in ride-hailing scenarios. The results are encouraging, especially when predictions are accurate, a case which is often encountered in practice when historical data is abundant. }
\end{abstract}

\section{Introduction}\label{sec1}

\rev{
Continuously varying optimization programs have appeared as a natural extension of time-invariant ones when the cost function, the constraints, or both, depend on a data streams and change continuously in time. This setting captures relevant problems in control, robotics, and signal processing, see, e.g.~\cite{DallAnese2020,Simonetto2020}. A typical example of such setting occurs when one want to determine the optimal operations of a smart grid with high renewable energy penetration, but they depend on constantly fluctuating supply~\cite{DallAnese2016}. Other typical examples range from video processing~\cite{Hamam2022} to robot control~\cite{Zavlanos2013}.

It becomes then important to study time-varying optimization problems of the form}
\begin{equation}~\label{eq:tv}
\min_{\x \in \reals^n} f(\x; \y(t)) + g(\x),\qquad t\geq 0,
\end{equation}
where $f: \reals^n \times \reals^d \to \reals$ is a smooth strongly convex function in $\x$ (for all $\y$), parametrized over a time-varying data stream $\y(t) \in \reals^d$ ($t$ represents the continuous time), and $g$ is a closed convex and proper function (such as the indicator function, or an $\ell_1$ regularization). \rev{Indeed, Problem~\eqref{eq:tv} models the need for optimizing for $\x$ when the cost function changes in time due to external data streams $\y(t)$, e.g., the fluctuating energy supply.  }


Solving Problem~\eqref{eq:tv} means to find and track the optimizer trajectory $\x^\star(t)$, as $\y(t)$ changes and it is revealed. In order to accomplish this, in the standard time-varying literature, one can sample Problem~\eqref{eq:tv} at discrete time instant $t_k$, $k=0,1,\ldots$, with $h:=t_{k}-t_{k-1}$ being the sampling time, and solve the sequence of time-invariant problems
\begin{equation}~\label{eq:tv:sampled}
\x^{\star}(t_k) = \arg\min_{\x \in \reals^n} f(\x; \y(t_k)) + g(\x), \qquad k\in \naturals,
\end{equation}
as they are revealed in time. Then one can set up online algorithms that find approximate $\x^{\star}(t_k)$'s, say $\x_k$'s, that eventually converge to\footnote{Somebody could rather say: ``track''.} the solution trajectory, within an error bound. The reader is referred to the surveys~\cite{DallAnese2020,Simonetto2020} for an ample treatment of the methods in both discrete and continuous time. One of the important aspects to keep in mind is that online algorithms are sought that are computationally frugal, so that one can approximate the solution of $\x^{\star}(t_k)$ within the sampling time $h$, and the key performance metric is how good the algorithms are with respect to an algorithm that has infinite computational time. 

A tacit assumption in the above methods is that one wants to converge to the solution trajectory generated by the evolution of the data stream $\y(t)$. However, this may be not the best course of action, since the data are often noisy and convergence to a noisy optimizer may be not advisable. A better angle is to ask whether we can set up algorithms to converge to a \emph{filtered} version of the trajectory. 

In this context, we re-interpret Problem~\eqref{eq:tv} as a stochastic problem, where we would like to find and track the filtered solution trajectory as
\begin{equation}~\label{eq:tv:expected}
\hat{\x}^{\star}(t) = \arg\min_{\x \in \reals^n} {\E}_{\y \in \mathcal{Y}(t)}[f(\x; \y)] + g(\x),
\end{equation}
where the expectation ${{\E}}_{\y \in \mathcal{Y}(t)}[\cdot]$ is with respect to the random variable $\y$, which is drawn from a time-varying probability distribution $\mathcal{Y}(t)$. 

Were the probability distribution time-invariant, Formulation~\eqref{eq:tv:expected} would be common in stochastic optimization. With our setting, we are considering instead a gradual distribution shift of an unknown distribution, which renders the formulation less common and more challenging. Recent papers have started to look into this formulation~\cite{Subhro2020,Cutler2021,Maity2022}; especially the third one is the closest to our goal, \rev{in which the authors adapt the standard prediction-correction algorithms to the stochastic setting by properly tuning the prediction and the correction step sizes and deriving convergence results. However, despite their careful analysis,} they do not consider a non-smooth component as our function $g(\x)$. 

In this paper, we look from a different angle and ask whether we can use a perturbed version of the optimality condition as a suitable dynamical model to do filtering. This will have the advantage to combine prediction and correction in novel and more performant ways. 
To fix the ideas on this novel angle, consider the unconstrained problem:
\begin{equation}
\x^{\star}(t) = \arg\min_{\x \in \reals^n} f(\x; \y(t)),
\end{equation}
with a strongly convex, doubly differentiable, and smooth function $f$ in $\x$. By perturbing the optimality condition $\nabla_{\x}f(\x^{\star}(t); \y(t)) = {\bf 0}_n$, we can derive the ordinary differential equation that describes how the optimizer evolves as~\cite{Polyak1987}:
\begin{eqnarray}\label{eq.dyn}
\frac{\mathrm{d}}{\mathrm{d} t}\x^{\star}(t) &=\ & - [\nabla_{\x\x} f(\x^{\star}(t); \y(t))]^{-1} \nabla_{\y\x} f(\x^{\star}(t); \y(t))\frac{\mathrm{d}}{\mathrm{d} t}\y(t) \nonumber\\&=:& F(\x^{\star}(t), \y(t), \dot{\y}(t)).
\end{eqnarray}
This is a nonlinear dynamical system. To derive a filter, we need to couple this model with a measurement equation, which tells us how far from the optimizer trajectory we are. We can thus use,
\begin{equation}\label{eq.meas}
\z(t) = \nabla_{\x} f(\x^{\star}(t); \y(t)),
\end{equation} 
as a measurement equation (i.e., it will be different from zero at any point but on the optimizer trajectory, where $\z(t) = {\bf 0}_n$). 

Armed with~\eqref{eq.dyn} and \eqref{eq.meas}, we could design a dynamical filter to reconstruct ${\x}^{\star}(t)$ based on a noisy data stream ${\y}(t)$.

\subsection{Contributions}

Starting from the continuous-time intuition from~\eqref{eq.dyn} and \eqref{eq.meas}, we develop discrete-time filters for unconstrained and constrained time-varying problems. In particular,

$\bullet$ We derive an extended Kalman filter for unconstrained and differentiable convex problems in the discrete-time setting starting from an algorithmic viewpoint;

$\bullet$ We generalize the filtering procedure for constrained and non-differentiable problems leveraging the nonlinear dynamical systems and nonlinear measurement equations coming from forward-backward algorithms and fixed-point residuals. We are then able to derive the optimal ``Kalman-style'' gain via both a worst-case approach and via dissipativity theory. We also present a possibly less conservative linear parameter-varying (LPV) methodology. 

$\bullet$ We showcase the benefit of the approach with respect to state-of-the-art prediction-correction methods (which our filters generalize), in numerical simulations stemming from a ride-hailing example with multiple companies in New York City.

\subsection{Related work}

Time-varying and stochastic optimization are vibrant research fields, and we do not plan to give an exhaustive account here. The reader is referred to~\cite{Fazlyab2016,DallAnese2020,Simonetto2020} for the first and \cite{Flaxman2005,Duchi2011,Schmidt2017,Duchi2018} for a sub-sample of the second, and the many references therein. 

\rev{The closest references to our current work in the time-varying domain are~\cite{Bastianello2020, Bastianello2023}. In~\cite{Bastianello2020} the authors develop a series of prediction-correction methods to find a track the optimizer trajectory. In particular, they build an extrapolation-based prediction, which we will use as well, and they give very general convergence bounds, which can then be leveraged for a variety of tracking algorithms. The main drawback of~\cite{Bastianello2020} is the fact that the authors optimize for a noisy decision variable and not for the more reasonable filtered version, as in Problem~\eqref{eq:tv:expected}. The works~\cite{Paper3,Qi2019} are similar to~\cite{Bastianello2020}. In \cite{Bastianello2023}, an internal-model principle is used to go beyond prediction-correction algorithms and derive controllers to track optimizers with zero asymptotical error. The approach of the paper is quite interesting, even though still looking at the noisy decision variables and not at their filtered version; moreover, one has to assume to have a model for the variation of the cost, which in principle we do not assume. Our approach can be seen as ``dual'' to a control design, since we design filters.    

In the stochastic setting and especially machine learning, the fact to have a distribution shift has led to the concept of continual learning, which is gaining traction, see~\cite{aljundi2019gradient, farajtabar2020orthogonal, lesort2020continual, lesort2021understanding}. All these works are focused on how to learn and re-learn machine learning losses via complex architectures. So, while very interesting and related, our methods differ in the basic assumption to have a convex optimization problem.} 

As such, \emph{Jointly} Stochastic \emph{and} time-varying convex optimization is a less studied area, and we have \rev{cited}~\cite{Subhro2020,Cutler2021,Maity2022} \rev{in the introduction}. \rev{In particular, the pre-print  of J.~Zhang and co-authors~\cite{Subhro2020}, and the published follow up~\cite{JZhang2022}, propose an online stochastic gradient algorithm that estimates the noise moments. This algorithm does not predict ahead but uses a filter to estimate the noise; moreover, no contraints nor non-smooth costs are considered. Finally, their analysis is at fixed time horizon $T$, requiring a restarting procedure for streaming data, which is not ideal in our setting. In~\cite{Cutler2021}, the authors propose an online proximal gradient algorithm to track a time-varying optimizer under noise drift. In this context, their setting is similar to ours. They look at no-prediction schemes (even thought with the addition of a linear smoothing term in the case of tracking of the optimal value) and they analyze convergence in expectation and high-probability. This particular paper is also close in techniques to inexact stochastic methods, such as~\cite{bastianello2021stochastic,Dixit2018}.  

The closest to our setting is then the work of~\cite{Maity2022} which further improves upon~\cite{Cutler2021} (at least in the differentiable setting) by proposing a prediction-correction method to find and track time-varying filtered optimizers. The authors carefully build a predictor based on a moving average method and tune the step sizes of prediction and correction as a function of the sampling period to obtain low asymptotical errors. They design the moving average by analyzing the trade-off between bias and variance. While this is very enticing, the overall methodology does not generalize immediately to non-smooth settings and one may want to decide on the step sizes and the sampling period independently, as we propose. }

The celebrated paper~\cite{Besbes2015} is also somewhat in the general area of interest, though their approach does not include prediction, it uses restart and is more directed at finding optimal regret rates for non-stationary objectives rather than noise. 

We will then build on techniques from dissipativity theory for analyzing and designing optimization algorithms. This is a recent and fertile area brought to fame by L. Lessard and collaborators' seminal paper~\cite{Lessard2016}, and now gaining momentum~\cite{Nelson2018,Colombino2019,Hassan-Moghaddam2021,Scherer2021,Michalowsky2021,Lessard2022}. \rev{The reader is also referred to the very recent works~\cite{syed2023bounds,notarnicola2023gradient}.} \rev{All the mentioned papers formulate the optimization algorithms as a nonlinear dynamical system, whose stability and performance need to be characterized via dissipativity analysis. This approach is often far from trivial.}
Our novel insight in this direction is to use the techniques to determine Kalman-style gains for optimization algorithms with errors. \rev{Our key step to ease the analysis is to be able to lock up most of the algorithm into two dissipative terms, and expose the gain between these, which we then optimize. } Finally, we employ LPV techniques from~\cite{Wu2005,Scherer2006,Massioni2020}, especially in the context of matrix sum-of-squares relaxations.  

\vspace*{2mm}
{\bf Organization.} 
The remaining of the paper is organized as follows. Section~\ref{sec:form} discusses formulation and main assumptions. We focus on the unconstrained case in Section~\ref{sec:ekf}, and on the general case in Section~\ref{sec:gen}. Section~\ref{sec:BMI} describes our gain design. We then conclude with some numerical simulations in Section~\ref{sec:num}. All proofs are given in the Appendix.

{\bf Notation.} Notation is wherever possible standard. For a differentiable function $f$, we define a step of the gradient method starting from a point $\x_k$ as $\x_{k+1} = [\OI - \alpha \nabla_{\x}f(\bullet)] \x_k \equiv \x_k - \alpha \nabla_{\x}f(\x_k)$, where $\alpha>0$ is the step size and $\OI$ is the identity operator. Then, $s$ steps are indicated as,
\begin{equation}\label{shn}
\x_{k+1} = \x_k^s = [\OI - \alpha \nabla_{\x}f(\bullet)]^{\circ s} \x_k.
\end{equation}
We let also $\x_k^0 = \x_k$ when needed. 

We further indicate with $\prox_{\alpha g}(\x)$ the proximal operator,
\begin{equation}
\prox_{\alpha g}(\x) = \arg\min_{\v \in \reals^n} \Big\{g(\v) + \frac{1}{2 \alpha}\|\v-\x\|^2\Big\}.
\end{equation}

Finally, spaces are indicated as $\reals, \naturals$, probability distributions are calligraphic, i.e., $\mathcal{Y}$, matrices and vectors are boldfaced, e.g., $\x \in \reals^{n}, \A \in \reals^{n\times m}$, operators are in sans-serif, e.g., $\OI, \OJ$, constants are in standard roman.  

\section{Problem formulation and assumptions}\label{sec:form}

Let us now consider the sequence of problems,
\begin{equation}~\label{eq:tv:sampled:unc}
\hat{\x}^{\star}(t_k) = \arg\min_{\x \in \reals^n} \E_{\y \in \mathcal{Y}(t_k)} [f(\x; {\y})] + g(\x), \qquad k\in \naturals,
\end{equation}
with a strongly convex, doubly differentiable and smooth cost function $f$ uniformly in $\x$ and a generic convex function $g$. Let us also introduce the shorthand notation $\hat{\x}^{\star}_{k} = \hat{\x}^{\star}(t_k)$, $\y_{k} = \y(t_k)$, and $\mathcal{Y}_k= \mathcal{Y}(t_k)$. Notice that, as done in stochastic optimization, the data point $\y_{k} = \y(t_k)$ is supposed to be a random vector drawn from the distribution $\mathcal{Y}_k= \mathcal{Y}(t_k)$.

First, let us derive a discrete-time dynamical system on how the optimizers evolve in time. Quite naturally, one could be attempted at discretizing~\eqref{eq.dyn}, but the presence of the inverse of the Hessian and the derivative of the data stream makes it quite cumbersome, especially if one has then to linearize it for an extended Kalman filter.

Instead we use a Bayesian approach, assuming that we have a (noisy) prior on how the data stream evolves, and we start from what we can be computed algorithmically. Let us denote with  $\OJ_{k+1}(\x): \reals^n \to \reals^n$ an approximation at time $t_{k}$ of ${\E}_{\y \in \mathcal{Y}_{k+1}}[\nabla_{\x} f(\x; \y)]$. Noisy prior on how the data evolves and how the gradient evolves can come from linear filters, or more sophisticated neural network models, or kernel models (see also~\cite{Rakhlin2013} for what they call \emph{predictable sequences}). For now think about one of the simplest model: $\OJ_{k+1}(\x) = 2 \nabla_{\x} f(\x; \y_{k}) - \nabla_{\x} f(\x; \y_{k-1})$ (obtained via an extrapolation technique~\cite{Bastianello2020}).

Then, we use an algorithmic view. At time $k$, we would like to solve
\begin{equation}~\label{eq:dyn:1}
\hat{\x}^{\star}_{k+1} = \arg\min_{\x \in \reals^n} \E_{\y \in \mathcal{Y}_{k+1}}[f(\x; \y)] + g(\x),
\end{equation}
yet it is not possible with data up to time $t_k$. So, we introduce the noisy dynamical system
\begin{equation}~\label{eq:dyn:2}
\hat{\x}^{\star}_{k+1} = [\prox_{\alpha g} (\OI-\alpha \OJ_{k+1}(\bullet))]^{\circ P} \hat{\x}^{\star}_{k} + \q_k =: \Phib_{k,g}(\hat{\x}^{\star}_{k})+ \q_k, 
\end{equation}
where $[\prox_{\alpha g} (\OI-\alpha \OJ_{k+1}(\bullet))]^{\circ P}$ means the application of the proximal gradient method of step size $\alpha>0$ for $P$ times, and $\q_k$ is the process error. The error comes both from a modelling error (truncating after $P$ iterations), and from the noisy predicted gradient $\OJ_{k+1}(\x)$. We will see how to characterize the error later, but for now it is useful to keep in mind that $\q_k$ is not-zero mean in general. If the gradient were exact and $P \to \infty$, then Equation~\eqref{eq:dyn:2} would solve~\eqref{eq:dyn:1} with no noise ($\q_k = 0$).

The ``pseudo''-dynamical system~\eqref{eq:dyn:2} will be our \emph{computationally affordable} nonlinear dynamical model. 

In par with~\eqref{eq:dyn:2}, we introduce a measurement equation,
\begin{eqnarray}~\label{eq:meas:2}
\z_{k+1} &=\ & -\hat{\x}^{\star}_{k} + [\prox_{\beta g}(\I-\beta \nabla_{\x} f(\bullet; {\y}_{k+1}))]^{\circ C} \hat{\x}^{\star}_{k} + \r_k \nonumber\\ &=:& \Psib_g(\hat{\x}^{\star}_{k},{\y}_{k+1})+ \r_k, 
\end{eqnarray}
where $C$ is the number proximal gradient steps, $\beta>0$ is the correction step size, and $\r_k$ is a noise term coming from the noisy character of $\y_{k+1}$ and it is in general not zero-mean. Again, $\z_{k+1} = {\bf 0}_n$ on the optimizer trajectory.

The right-hand side of~\eqref{eq:meas:2} represents the fixed-point residual of our $C$-steps proximal gradient method, that we use to compute the measurement. 

\subsection{Properties, requirements of the gradient approximators}

Before going further, it is useful to understand a bit better the properties of the gradient approximations $\OJ_{k+1}(\x)$ and $\nabla_{\x}f(\x; \y_{k+1})$. Both are attempting at approximating $\E_{\y \in \mathcal{Y}_{k+1}}[\nabla_{\x} f(\x; \y)]$, but there are a few differences.

The first is based on data available up to time $t_k$ and it is in general a biased estimator. This is not a problem per se, since even in deterministic prediction-correction methods, the predicted gradient is in general a biased prediction. 

\begin{example}[Recurring stochastic example]\label{ex:2}
Consider the case in which $\nabla_{\x}f(\x;\y)$ is linear in $\y$ (e.g., for linear models and Least-Squares estimators, for example when $f(\x;\y):= \frac{1}{2}\|\x-\A\y\|^2$ for a given matrix $\A$), and more generally the case in which:
$$
f(\x;\y) = f'(\x) + \y^\transp \A \x,
$$
with $f'(\x)$ strongly convex and smooth. In this case, let $\|\nabla_{\y\x}f(\x;\y)\| = \|\A\| \leq C_0$. 

Further, suppose $\y(t)$ is generated by a nominal (doubly differentiable in $t$) trajectory to which we add a Gaussian zero-mean noise at each sampling time $t_k$. And in particular, set $\y(t_k) = \bar{\y}(t_k) + \e_k$, with $\e_k \sim \mathcal{N}({\bf 0}, \Sigmab_k)$ for a given time-varying covariance matrix $\Sigmab_k$. We also know that $\E[\|\e_k\|] \leq \sqrt{\trace(\Sigmab_k)}$, and we set $\sqrt{\trace(\Sigmab_k)} \leq \Sigma$.  

\rev{We consider now three different prediction strategies. In particular, we consider using the gradient of the current function to predict the gradient of the future function, we consider an extrapolation prediction with two points and three points. For the extrapolation formulas, we refer to~\cite{Bastianello2020}. We indicate with the superscript $(1)$, $(2)$, $(3)$ the three strategies.

We then have,  
\begin{eqnarray*}
\OJ_{k+1}^{(1)}(\x) &=& \nabla_{\x} f(\x; \y_{k}), \qquad 
\OJ_{k+1}^{(2)}(\x) = 2\nabla_{\x} f(\x; \y_{k}) - \nabla_{\x} f(\x; \y_{k-1}) = \nabla_{\x} f(\x; 2\y_{k} - \y_{k-1}), \\
\OJ_{k+1}^{(3)}(\x) &=& 3\nabla_{\x} f(\x; \y_{k}) - 3\nabla_{\x} f(\x; \y_{k-1}) + \nabla_{\x} f(\x; \y_{k-2})= \nabla_{\x} f(\x; 3\y_{k} - 3\y_{k-1}+\y_{k-2}).
\end{eqnarray*}
}

Consider a nominal trajectory $\bar{\y}(t)$ for which we assume the bounds:
\begin{equation}\label{bounds2}
\max\{\left\|\nabla_{t} \bar{\y}(t)\right\|, \left\|\nabla_{tt} \bar{\y}(t)\right\|, \left\|\nabla_{ttt} \bar{\y}(t)\right\|\} \leq C, \qquad \forall \x, t,
\end{equation}
\rev{where the need as many derivatives as number of points in the predictors. }
Then, in the Appendix we show that: \rev{
\begin{eqnarray*}
\E[\|\OJ_{k+1}^{(1)}(\hat{\x}^{\star}_{k+1}) - \E_{\y \in \mathcal{Y}_{k+1}}[\nabla_{\x} f(\hat{\x}^{\star}_{k+1}; \y)]\|] &=& \E_{\w \in \mathcal{Y}_{k}}[\|\nabla_{\x} f(\hat{\x}^{\star}_{k+1}; \w)  -\E_{\y \in \mathcal{Y}_{k+1}}[\nabla_{\x} f(\hat{\x}^{\star}_{k+1}; \y)]\|]  \\ &\leq& C_0 h + C_0 \Sigma, \\
\E[\|\OJ_{k+1}^{(2)}(\hat{\x}^{\star}_{k+1}) - \E_{\y \in \mathcal{Y}_{k+1}}[\nabla_{\x} f(\hat{\x}^{\star}_{k+1}; \y)]\|] &=& \E_{\w \in \mathcal{Y}_{k}, \z \in \mathcal{Y}_{k-1}}[\|\nabla_{\x} f(\hat{\x}^{\star}_{k+1}; 2\w - \z)  -\\&&-\E_{\y \in \mathcal{Y}_{k+1}}[\nabla_{\x} f(\hat{\x}^{\star}_{k+1}; \y)]\|]  \\ &\leq& C_0 C h^2 + 3 C_0 \Sigma, \\
\E[\|\OJ_{k+1}^{(3)}(\hat{\x}^{\star}_{k+1}) - \E_{\y \in \mathcal{Y}_{k+1}}[\nabla_{\x} f(\hat{\x}^{\star}_{k+1}; \y)]\|] &=& \E_{\w \in \mathcal{Y}_{k}, \z \in \mathcal{Y}_{k-1}, \q \in \mathcal{Y}_{k-2}}[\|\nabla_{\x} f(\hat{\x}^{\star}_{k+1}; 3\w - 3\z + \q)  -\\&&-\E_{\y \in \mathcal{Y}_{k+1}}[\nabla_{\x} f(\hat{\x}^{\star}_{k+1}; \y)]\|]  \\ &\leq& C_0 C h^3 + 7 C_0 \Sigma,
\end{eqnarray*} 
}
and
\begin{equation*}
\E_{\y\in\mathcal{Y}_{k+1}} [\|\nabla_{\x}f(\x; \y) - \E_{\y \in \mathcal{Y}_{k+1}}[\nabla_{\x} f(\x; \y)] \|] \leq  C_0 \Sigma.
\end{equation*}

\rev{
In the error bounds for the predictions, one can immediately see the trade-off between a bias term, which depends on powers of the sampling period, and a variance $O(\Sigma)$ term. Higher-order (more-memory) predictors have lower bias (since $h<1$), but higher variance. 
}
\qed
\end{example}

\smallskip
\smallskip
\smallskip

The second estimator, meaning using $\nabla_{\x}f(\x; \y_{k+1})$ instead of $\E_{\y \in \mathcal{Y}_{k+1}}[\nabla_{\x} f(\x; \y)]$, is evaluated with data coming at time $t_{k+1}$ and it is unbiased~\cite{Flaxman2005}.

For the two approximations we require the following.

\begin{assumption}\label{as.1}
Let the cost function $f(\x;\y)$ be $\mu$-strongly convex and $L$-smooth in $\x$ uniformly in $\y$ (i.e., for all $\y$). The chosen gradient predictor $\OJ_{k+1}(\x)$ is then $\mu$-strongly monotone and $L$-Lipschitz in $\x$ for all $k$'s. 
\end{assumption}

\begin{assumption}\label{as.2}
The noise processes and gradient prediction errors are bounded as follows:
\begin{eqnarray*}
(a) & \E [\|\OJ_{k+1}(\hat{\x}^{\star}_{k+1}) - \E_{\y \in \mathcal{Y}_{k+1}}[\nabla_{\x} f(\hat{\x}^{\star}_{k+1}; \y)]\|] \leq \tau, \\
(b) & \E_{\y\in\mathcal{Y}_{k+1}} [\|\nabla_{\x}f(\x; \y) - \E_{\y \in \mathcal{Y}_{k+1}}[\nabla_{\x} f(\x; \y)] \|] \leq \sigma, 
\end{eqnarray*}
for finite scalars $\tau$ and $\sigma$, for all $k \in \naturals$, and (b) for all $\x \in \reals^n$.
\end{assumption} 

Assumption~\ref{as.1} is often required for time-varying optimization~\cite{Simonetto2020}. The assumption on the predictor is also reasonable, for instance is verified in Example~\ref{ex:2}, and with Taylor-based and extrapolation-based predictions~\cite{Bastianello2020}. 

For Assumption~\ref{as.2}: Property $(a)$ is in par with some deterministic and stochastic assumptions appeared in past years. For example, it can be seen in parallel with the quality of the hint or predictable sequences in~\cite{Rakhlin2013,Jadbabaie2015}. In Example~\ref{ex:2} \rev{with predictor $(2)$}, Property $(a)$ is verified with $\tau = C_0 C h^2 + 3 C_0 \Sigma$. Stochastic versions of Property $(a)$ have appeared, e.g., in~\cite{Maity2022}, with example-based constructions for determining a suitable $\OJ_{k+1}$.   
Property $(b)$ is also commonly asked in stochastic frameworks~\cite{Subhro2020}. Usually, one asks that $\E_{\y\in\mathcal{Y}_{k+1}} [\|\nabla_{\x}f(\x; \y) - \E_{\y \in \mathcal{Y}_{k+1}}[\nabla_{\x} f(\x; \y)] \|^2] \leq \sigma^2$, but due to the convexity of $(\cdot)^2$ and Jensen's inequality, Property $(b)$ is implied by the squared one (in fact $(\E[\|\cdot\|])^2 \leq \E[\|\cdot\|^2]$). For us Property $(b)$ can be time-varying. In Example~\ref{ex:2}, Property $(b)$ is verified with $\sigma = C_0 \Sigma$. Finally, Property $(b)$ can be tightened to be valid only on the algorithm iterates $(\x_k)_{k\in\mathbb{N}}$, and interpret it as gradient noise with a small theoretical effort~\cite{Cutler2021}.

\section{An extended Kalman filter}\label{sec:ekf}

We are now ready to derive a filter to track the filtered optimizer trajectory $\hat{\x}^{\star}(t_k)$. Since the \rev{dynamical} system and the measurement equations are nonlinear, we will use an extended Kalman filter. For it, we require that both $\nabla_{\x}f(\x;\y)$ and $\OJ_{k+1}(\x)$ are differentiable with respect to $\x$, and we will require knowledge of the covariance of the noise processes $\q_k, \r_k$ at all time instances. We will also assume that $g \equiv 0$, to be able to differentiate, which puts us in an unconstrained differentiable problem setting, where we use gradient (the proximal operator is the identity in this case). 

Recall the notation $\x_{k}^{s} = [\OI - \alpha \nabla_{\x} f(\bullet, {\y}_{k+1})]^{\circ s}\x_{k}$ defined in~\eqref{shn}, indicating the effect of $s$ steps of the gradient method, or an approximate gradient if $f(\bullet, {\y}_{k+1})$ is substituted with $\OJ_{k+1}(\bullet)$. Define the derivative quantities (we drop the mention to $g$, since it does not exist in this case),
\begin{eqnarray}
\hspace*{-.25cm}\F_k &=&  \nabla_{\x} \Phib_k(\x_{k}) = \prod_{p=1}^P\left(\I_n-\alpha \nabla_{\x}\OJ_{k+1}(\x_{k}^{P-p})\right) \label{eq.F}\\
\H_{k+1} &=&  \nabla_{\x} \Psib(\x_{k+1|k}, \y_{k+1}) =  \I_n - \prod_{c=1}^C\left(\I_n-\beta \nabla_{\x\x}f(\x_{k+1|k}^{C-c}; \y_{k+1})\right). \label{eq.H}
\end{eqnarray}
We also let $\R_k \in \reals^{n\times n}$ and $\Q_k\in \reals^{n\times n}$ be the covariance matrices of the noise processes $\r_k$ and $\q_k$, respectively.

\begin{algorithm}
\caption{An extended Kalman filter (TV-EKF)}\label{alg:cap}
\begin{algorithmic}[1]
\Require Initialize: ${\x}_{1|0} = {\bf 0}$, $\P_{1|0} = \I_n$. Number of prediction and correction steps $P, C$, sampling time $h$, step sizes $\alpha, \beta$, covariance matrices $\R_k, \Q_k$ for all $k$, as well as prediction strategy $\OJ$.  
\Ensure A sequence $({\x}_k)_{k\in \mathbb{N}}$
\For{$k\in \naturals, k\geq1$}
\State Receive $\y_{k}$
\State Compute $\H_{k}$ as in Eq.~\eqref{eq.H}. 
\State Correction step: 
\begin{eqnarray*}
\K_k &=&  \P_{k|k-1}\H_{k}(\H_{k} \P_{k|k-1} \H_{k}^\transp + \R_k)^{-1}\\
\x_{k} &=& \x_{k|k-1} + \K_k (\Psib(\x_{k|k-1}, {\y}_{k})) \\
\P_{k} &=& (\I - \K_k \H_{k}) \P_{k|k-1}
\end{eqnarray*}
\State Compute $\F_k$ as in Eq.~\eqref{eq.F}.
\State Prediction step: 
\begin{equation*}
\x_{k+1|k} =  \Phib_k(\x_k), \qquad \P_{k+1|k} = \F_k \P_k \F_k^\transp + \Q_k
\end{equation*}
\EndFor
\end{algorithmic}
\end{algorithm}

With this in place, the extended Kalman filter (TV-EKF) represented in Algorithm~\ref{alg:cap} is able to filter and track the optimizer trajectory $\hat{\x}^{\star}(t)$. We notice that we have presented the filter with correction first, to highlight the standard workflow within a sampling period. We notice also that the filter can be extended to include a filtering process for the data stream $\y(t)$, if a dynamical model for it is available. 

Algorithmically, the presented TV-EKF requires several computations. In the correction step, it comprises the computation of the Hessian of $f$ at various points for determining $\H_k$ and taking a matrix inverse for determining $\K_k$. The update for $\x_k$ requires also $C$ gradient steps. In the prediction pass, the filter includes a process to determine any prediction $\OJ_{k+1}$, computing its derivatives, and $P$ gradient steps. Considering a $n$-dimensional state, and letting $\mathsf{g}, \mathsf{h}, \mathsf{j}, \mathsf{dj}$ be the computational effort to determine the gradient, Hessian, predicted gradient, and its derivatives, then the overall computational complexity of TV-EKF is $O(C(\mathsf{h} + \mathsf{g}) + P(\mathsf{dj} + \mathsf{g}) + n^3)$. 

We will study the empirical performance of TV-EKF in Section~\ref{sec:num}, but we close here with some remarks. 

\begin{remark}[Prediction-Correction methods]\label{rem:pc} We can see how $\x_k$ is updated as 
\begin{multline}
\x_{k} = \Phib_k(\x_{k-1}) + \K_{k} (\Psib(\Phib_k(\x_{k-1}); \y_{k})) =  
 (\I_n - \K_k)\underbrace{[\OI-\alpha \OJ_{k}(\bullet)]^{\circ P}\x_{k-1}}_{\textrm{prediction}}
  + \\+ \K_{k}\underbrace{[\OI-\beta \nabla_{\x} f(\bullet; {\y}_{k})]^{\circ C} \circ[\OI-\alpha  \OJ_{k}(\bullet)]^{\circ P}\x_{k-1}}_{\textrm{correcting the prediction}},
\end{multline}
and if we let $\K_k \equiv \I_n$, then we obtain back the standard prediction-correction methods. \rev{In the latter case, we will see its asymptotical error as a special case of Theorem~\ref{prop.2}.}
\end{remark}
 
\begin{remark}
The TV-EKF algorithms that is presented here could be extended to equality-constrained optimization problems, once formulated as saddle-points, but we leave this for future endeavors. 
\end{remark}

The TV-EKF algorithm that we have presented offers several advantages, and above all the ease of implementation. However, from an optimization perspective, it is lacking in good convergence guarantees\footnote{Besides the fact that the prediction and correction steps represent contractive operators for $\alpha<2\mu/L^2, \beta<2/L$.}, and from a noise perspective, we do not have a good intuition or recipe on how to set the covariances $\R_k, \Q_k$ in meaningful ways. 

This, combined with the fact that TV-EKF will not work for non-smooth costs, pushes us to look beyond to a more general setting.  However, before moving on, we give some intuition on why the TV-EKF does perform well empirically in the numerical settings that are presented in Section~\ref{sec:num}.

\begin{proposition}[Equivalence to a damped Newton's step]\label{pr.kalman}
When the noise on the measurement is negligible, i.e., $\R_k \approx 0$, and we take only one step of correction $C = 1$, then TV-EKF is a damped Newton's method, with update 
\begin{equation*}
\x_{k} = \x_{k|k-1} - \beta[\nabla_{\x\x}f(\x_{k|k-1}; {\y}_{k})]^{-1} \nabla_{\x} f(\x_{{k|k-1}}; {\y}_{k}). 
\end{equation*}\qed
\end{proposition}

Proposition~\ref{pr.kalman} implies that TV-EKF includes second-order information and could be thought of as a stochastic quasi-Newton method. In this sense, if the noise covariances are well estimated, or small, then TV-EKF is expected to do better than standard prediction-correction methods that only use first-order information. This is not surprising, but the connection Kalman-Newton in optimization is interesting. \rev{With this in place, we now expect that the Extended Kalman filter approach can offer better performances then a standard approach. So, it is natural to see whether one could use a filtering approach also in more general optimization settings, as we explore next. }

\section{The general case}\label{sec:gen}

The standard extended Kalman filter can be easily derived when the cost is differentiable. We generalize now our filter design to the case in which one has also the term $g(\x)$ in the cost, modeling constraints and non-smooth regularizations. 

Reconsider the dynamical system~\eqref{eq:dyn:2} in par with the measurement equation~\eqref{eq:meas:2}. Under Assumption~\ref{as.1}, from standard operator theory, we know that if $\alpha$ and $\beta$ are chosen small enough, and in particular\footnote{This is due to the fact that $\OJ_{k}$ is $\mu$-strongly monotone and $L$-Lipschitz, but not a gradient per se. Sharper conditions can be derived if $\OJ_{k}$ would be a gradient, like in the correction case, for which we can choose $\alpha<2/L$. } $\alpha<2\mu/L^2, \beta<2/L$, both the prediction and the correction represent contraction operators, which converge to their respectives unique fixed points. In particular, their contraction factors are~\cite{Nesterov2004}: 
\begin{equation}\label{contra.factors}
\rhop = \sqrt{1-2\alpha \mu+\alpha^2L^2}, \qquad \rhoc = \max\{|1-\beta\mu|,|1-\beta L|\}.
\end{equation}
for prediction and correction operators, respectively. 

\subsection{A static gain filter}

With our dynamical model~\eqref{eq:dyn:2} and measurement equation~\eqref{eq:meas:2}, we are now ready to build our filter. Here, since the model and the measurements are non-differentiable equations, \rev{an Extended Kalman filter approach cannot be readily applied, since the matrices $\F_k$ and $\H_k$ in \eqref{eq.F} and \eqref{eq.H} cannot be computed}. \rev{Hence,} we will focus on static gain filters, that can be computed off-line, before running the time-varying algorithm. We let $\Psib_{g}'(\x, \y) = \Psib_{g}(\x,\y) + \x$. 

Therefore, we start by considering the update equation
\begin{equation}
\x_{k+1} = \Phib_{k,g}(\x_k) - \K (\Phib_{k,g}(\x_k)- \Psib_{g}'(\Phib_{k,g}(\x_k), \y_{k+1} ) )=
(\I_n-\K)\Phib_{k,g}(\x_k) + \K \Psib'_{g}(\Phib_{k,g}(\x_k) , \y_{k+1} ),\qquad k \in \naturals
\end{equation}
consisting of running a prediction $\Phib_{k,g}(\x_k)$ and then correcting it via the correction $\Psib'_{g}(\Phib_{k,g}(\x_k), \y_{k+1} )$. As mentioned in Remark~\ref{rem:pc}, by setting $\K\equiv\I_n$, we obtain back the standard prediction-correction methods. 

Algorithm~\ref{alg:contr} makes the update explicit along with all the involved computations, for a generic choice of $\OJ$ and $\K$. As we can see the computational complexity here does not involve matrix inversions, but it adds proximal mapping computations. If the proximal step is easy to perform compared to the other computations, \rev{for example when the cost is prox-friendly and the proximal step can be carried out in closed-form}, then the complexity is $O(C(\mathsf{h} + \mathsf{g}) + P(\mathsf{dj} + \mathsf{g}))$, which is better than TV-EKF, as expected.

\begin{algorithm}
\caption{A static contractive filter (TV-CONTRACT)}\label{alg:contr}
\begin{algorithmic}[1]
\Require Initialize: ${\x}_{1|0} = {\bf 0}$. Number of prediction and correction steps $P, C$, sampling time $h$, step sizes $\alpha, \beta$, prediction strategy $\OJ$, as well as a filter gain $\K$.  
\Ensure A sequence $({\x}_k)_{k\in \mathbb{N}}$
\For{$k\in \naturals, k\geq1$}
\State Receive $\y_{k}$ 
\State Correction step: $\x_k = (\I_n-\K)\x_{k|k-1} + \K \Psib_{g}'(\x_{k|k-1} , \y_{k})$
\State Compute $\OJ_{k+1}(\x_k)$ 
\State Prediction step: Compute $\x_{k+1|k} = \Phib_{k,g}(\x_k)$
\EndFor
\end{algorithmic}
\end{algorithm}

As mentioned before, in this paper, we are interested in designing $\K$ in such a way to reduce the tracking error of the sequence $\{\x_k\}_{k \in \naturals}$, and in particular: which $\K$ would deliver the smallest $\limsup_{k\to \infty}\E[\|\x_k - \hat{\x}^\star_k\|]$? 

\subsection{Scalar, worst-case convergence results}\label{sec:wcg}

We start by analyzing the easier case of determining the best \emph{scalar} gain $\chi \in [0,1]$, for the update,
\begin{equation}\label{scalar}
\x_{k+1} =
(1-\chi)\Phib_{k,g}(\x_k) + \chi \Psib_{g}'(\Phib_{k,g}(\x_k) , \y_{k+1} ), \quad k \in \naturals.
\end{equation}
While this is restrictive in practice, it will give us some intuition on the general problem. Also, by restricting $\chi$ in $[0,1]$, we are considering all convex combinations of prediction and correction phases. The latter is important if $g$ represents a feasible set and we want the sequence $\{\x_{k}\}_{k \in \naturals}$ to be feasible for every $k$. 
We have the following theorem. 

\begin{theorem}\label{prop.2}
Let Assumptions~\ref{as.1}-\ref{as.2} hold. Assume furthermore that the optimizer trajectory is bounded as,
\begin{equation}\label{as.3}
\|\hat{\x}^\star_{k+1} - \hat{\x}^\star_{k}\| \leq \Delta < \infty, \quad \forall k \in \naturals.
\end{equation}
Choose $\alpha<2\mu/L^2, \beta<2/L$. Consider Algorithm~\ref{alg:contr} with the selection of $\K = \chi \I_n$, $\chi \in [0,1]$, leading to the update~\eqref{scalar}, and its sequence $\{\x_k\}_{k \in \naturals}$.
Define functions $\zeta$ and $\xi$ as:
\begin{eqnarray*}
&&\zeta_{\ell, \rho} = \{1 \textrm{ if } \ell = 0; \quad \rho^\ell \textrm{ otherwise }\}, \\ 
&&\xi_{\ell, \rho} = \{0 \textrm{ if } \ell = 0; \quad 1+\rho^\ell \textrm{ otherwise }\}. 
\end{eqnarray*}
Recall the contraction parameters~\eqref{contra.factors} and choose the number of prediction and correction steps $P$ and $C$ such that $\zetac \zetap < 1$.
Then, by calling $\varrho_\chi = (1-\chi) \zetap + \chi \zetac \zetap$, the asymptotic error is upper bounded as,
\begin{equation}\label{ate.w}
\limsup_{k \to \infty}\, \E[\|{\x}_{k} - \hat{\x}^\star_{k}\|] = \frac{1}{1-\varrho_\chi}\Big[(1-\chi)\left[\zetap \Delta + \xip \tau_\mu \right] + \chi \left[\zetac[\zetap \Delta + \xip \tau_\mu] + \sigmac\right]\Big],
\end{equation}
with $\tau_\mu = \tau/\mu$, and $\sigma_c = \beta \sigma/(1-\rhoc)$

Finally, under the setting of Example~\ref{ex:2}, \rev{predictor $(2)$}, $\Delta = C_0 h/\mu$, $\tau_{\mu} = (C_0 Ch^2 + 3C_0\Sigma)/\mu$ and $\sigma = C_0 \Sigma$. \qed
\end{theorem}

Theorem~\ref{prop.2} captures the asymptotic tracking error of the proposed TV-CONTRACT algorithm, when $\K = \chi\I_n$. The requirement~\eqref{as.3} is standard, assuring that the trajectory is regular enough to be tracked, see~\cite{Simonetto2020}. The condition $\zetac \zetap < 1$ is verified, whenever $P+C\geq1$, since $\rhop, \rhoc \in [0,1)$ with the choice of $\alpha, \beta$.  

For $\chi = 1$, Theorem~\ref{prop.2} extends \cite[Proposition~1]{Bastianello2020} to stochastic settings and \cite[Theorem~2.7]{Maity2022} to $f+g$ settings with multiple prediction and correction steps. \rev{In particular, we see how if the number of steps $P$ and $C$ increases, then the expected worst-case error decreases.}

The question we have for filter design is {\bf how to tune $\chi$ to lower the asymptotical error, given all the rest fixed?} 

\subsection{Tuning $\chi$}

The filter design problem can be now formulated as,
\begin{equation}\label{p.w}
\min_{\chi \in [0,1]}\, \eqref{ate.w}
\end{equation}
Problem~\eqref{p.w} is a linear-fractional programming, that can be solved by transforming it into a linear program, once all the coefficients are fixed. We do not give the details of this, since easily found in standard books~\cite[Ch. 4.3.2]{Boyd2004a}. Nonetheless we report an interesting fact on the nature of the solution. 

\begin{proposition}\label{lemma.3}
The solution of the tuning problem~\eqref{p.w} is either $\chi^{\star} = 1$, $\chi^{\star} = 0$, or in a special case, any $\chi \in [0,1]$. \qed
\end{proposition}

What Proposition~\ref{lemma.3} says is that from a worst-case perspective (i.e., from an asymptotic tracking error) we are better off to either just do predictions, or just do prediction-correction (and in a very special case, we can take any choice). The choice is made a priori, from the size of the prediction or correction errors (see the proof for exact conditions). 

This is hardly satisfactory. We see next how to extend the above to a generic matrix gain $\K$, via dissipativity theory, which has a richer behavior.

\section{Dissipativity theory and filter design}\label{sec:BMI}

We move now to the general case of designing a matrix $\K$ in an optimal fashion. We will be using recent tools from dissipativity theory applied to optimization algorithm design, and we were particularly inspired by~\cite{Colombino2019,Hassan-Moghaddam2021}.

\subsection{Filter design}

To start our filter design, we need to recast Algorithm~\ref{alg:contr} as a block diagram, where the optimization algorithmic updates are interpreted as nonlinear blocks and modeled as quadratic constraints. 
Consider then noisy algorithmic update,
\begin{equation}\label{algo.update}
\left\{\begin{array}{ll}
\w_{k+1} = \bar{\w}_{k+1} + \Q \q_{k+1}, & \bar{\w}_{k+1} = T_{k+1}^1(\x_k) \\
\uu_{k+1} = \bar{\uu}_{k+1} + \R \r_{k+1}, & \bar{\uu}_{k+1} = T_{k+1}^2(\bar{\w}_{k+1}) \\
\x_{k+1} = (\I-\K) {\w}_{k+1} + \K \uu_{k+1} &
\end{array}\right.
\end{equation}
where $\q_{k+1}$ and $\r_{k+1}$ are noise terms, whose expected norm is bounded, as we will see shortly, and $\Q, \R$ are tuning matrices that can model the relative amount of error or correlations.

\begin{figure}
\centering
\includegraphics[width=0.6\textwidth]{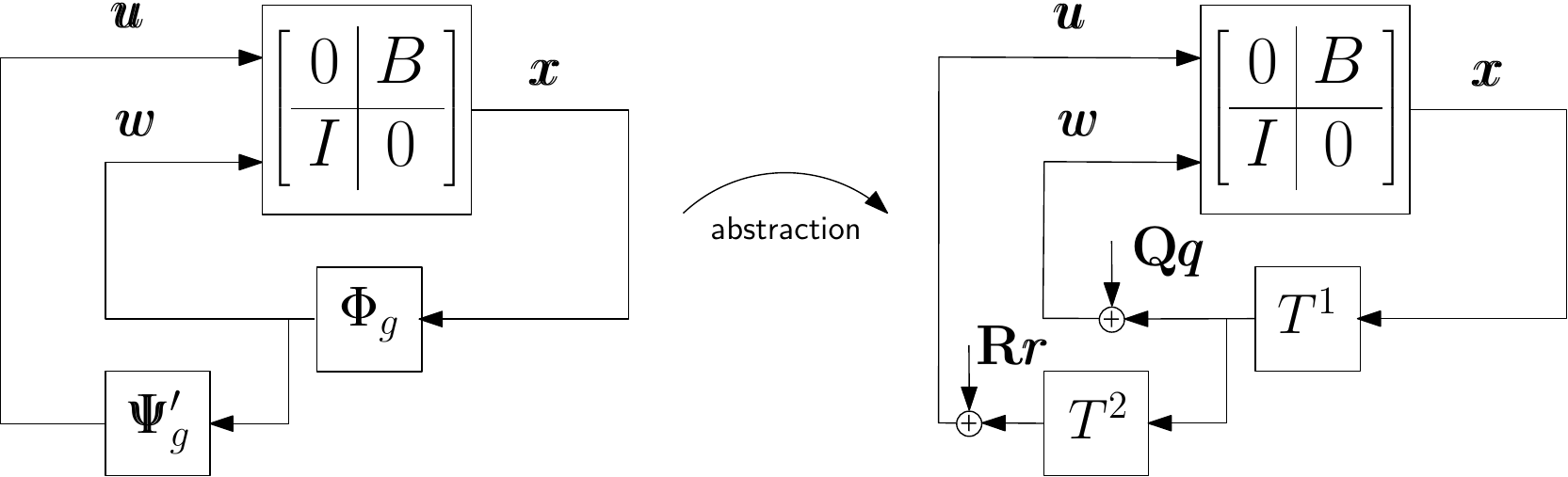}
\caption{The algorithmic choice~\eqref{algo.update} is an abstraction of Algorithm~\ref{alg:contr}. The matrix block indicates the algorithmic update.}
\label{fig.depict}
\end{figure}

The algorithmic choice~\eqref{algo.update} is an abstraction of Algorithm~\ref{alg:contr}, as we can see in Figure~\ref{fig.depict}, where $T^1_{k+1}$ and $T^2_{k+1}$ are the ideal operators representing the ideal prediction and correction steps, respectively, and both with fixed point $\hat{\x}^\star_{k+1} = \bar{\w}_{k+1}^\star = \bar{\uu}_{k+1}^\star$. The fact that, technically, one should consider $\bar{\uu}_{k+1} = T_{k+1}^2({\w}_{k+1})$, and the noise terms $\q_{k+1}$ and $\r_{k+1}$ are in fact correlated, since $\r_{k+1}$ should depend on a noisy prediction, can be ignored here, since we will only look at worst-case performance guarantees, from bounded errors to bounded output. 

Under Assumptions~\ref{as.1} and \ref{as.2},  with the same notation of Theorem~\ref{prop.2} and from the proof of \cite[Proposition 1]{Bastianello2020} with $\E[\|\tau_k\|] \leq \tau_{\mu}$, we can bound,
\begin{equation}
	\|{\w}_{k+1} - \bar{\w}^\star_{k+1}\| = \|{\x}_{k+1|k} - \hat{\x}^\star_{k+1}\|  \leq  \zetap\|{{\x}_{k} - \hat{\x}^\star_{k+1}}\| + \xip \tau_k.
\end{equation}
On the other hand, 
\begin{equation}
	\|{\w}_{k+1} - \bar{\w}^\star_{k+1}\| = \|\bar{{\w}}_{k+1} - \bar{\w}^\star_{k+1} + \Q \q_{k+1}\|  \leq \zetap\|{{\x}_{k} - \hat{\x}^\star_{k+1}}\| + \|\Q \q_{k+1}\|,
\end{equation}
so we can look at bounded errors as $\|\Q \q_{k+1}\| \leq \xip \tau_k$ and $\E[\|\Q \q_{k+1}\|] \leq \xip \tau_{\mu}$.

Similarly for the error term $\R \r_k$, we can look only for bounded errors as $\E[\|\R \r_k\|] \leq \zetac \xip \tau_{\mu} + \sigmac$. 

In particular, we consider $\E[\|\q_k\|]\leq 1/\sqrt{2}$, $\E[\|\r_k\|]\leq 1/\sqrt{2}$ and model the matrices $\Q, \R$ to have their largest singular value at $\sqrt{2}(\xip \tau_{\mu})$ and $\sqrt{2}(\zetac \xip \tau_{\mu} + \sigmac)$, respectively. 

For the sake of ease of notation, in~\eqref{algo.update}, we define matrices $\B, \B_e$, and the nominal input and error signal $\bar{\z}, \e$ as,
\begin{eqnarray}\label{algo.update}
&&\B = [(\I_n-\K), \K], \quad \B_e = [(\I_n-\K)\Q, \K\R],\quad\\ && \bar{\z} = [\bar{\w}^\transp, \bar{\uu}^\transp]^\transp, \quad {\e} = [{\q}^\transp, {\r}^\transp]^\transp, \nonumber \\
\label{abstraction}
&&\x_{k+1} = \B \bar{\z}_{k+1} + \B_e \e_{k+1},
\end{eqnarray}
and $\E[\|\e_k\|] \leq 1$. We also introduce point-wise-in-time quadratic constraints for $T^1$ and $T^2$ as,
\begin{equation}
\left[\!\begin{array}{c}  \x_{k} - \hat{\x}^\star_{k+1} \\\bar{\w}_{k+1} - \bar{\w}^\star_{k+1}\end{array}\!\right]^\transp \!\left[\!\begin{array}{cc}  \omega_1^2\I_n\!\! & 0\\ 0& \!\!-\I_n \end{array}\!\right]\!\left[\!\begin{array}{c}  \x_{k} - \hat{\x}^\star_{k+1} \\\bar{\w}_{k+1} - \bar{\w}^\star_{k+1}\end{array}\!\right] \geq 0,
\end{equation}
\begin{equation}
\left[\!\begin{array}{c}  \x_{k} - \hat{\x}^\star_{k+1} \\\bar{\uu}_{k+1} - \bar{\uu}^\star_{k+1}\end{array}\!\right]^\transp \!\left[\!\begin{array}{cc}  \omega_1^2\omega_2^2\I_n\!\! & 0\\ 0& \!\!-\I_n \end{array}\!\right]\!\left[\!\begin{array}{c}  \x_{k} - \hat{\x}^\star_{k+1} \\\bar{\uu}_{k+1} - \bar{\uu}^\star_{k+1}\end{array}\!\right] \geq 0,
\end{equation}
with $\omega_1 = \zetap$ and $\omega_2 = \zetac \zetap$, which are due to the contracting properties of $T^1$ and $T^2$, respectively. 
With this in place, convergence and asymptotic tracking error performance can be formulated as follows.

\begin{theorem}\label{th.3}
Consider Algorithm~\ref{alg:contr} and its abstraction~\eqref{abstraction}, to find and track the filtered optimizer trajectory $\hat{\x}^{\star}(t)$ of the time-varying stochastic optimization problem $\min_{\x\in\reals^n} \E[f(\x;\y(t))]+g(\x)$. Assume that the optimizer trajectory varies in a bounded way as $\|\hat{\x}^{\star}(t_{k+1})- \hat{\x}^{\star}(t_k)\| \leq \Delta$, for all $k\in \naturals$ and $\Delta<\infty$. Let Assumptions~\ref{as.1}-\ref{as.2} hold as well.  

Introduce matrix $\W\in\reals^{n\times n}$ and $\X\in\reals^{n\times n}, \X\succ 0$, scalars $\lambda_1\geq0, \lambda_2\geq0$, in addition to supporting scalars $\gamma_1, \gamma_2$, and consider the tuning matrix $\K$. For any fixed scalar $\rho\in(0,1)$, by solving the convex problem\footnote{For readability, we indicate in {\color{blue}blue} the decision variables.}
\begin{subequations}\label{condit.3}
\begin{eqnarray}
\hskip-1cm\minimize_{\scriptsize\begin{array}{c}\Xcb\succ 0, \lambdacb_1\geq0, \lambdacb_2\geq0,\\ \Wcb, \gammacb_1^2, \gammacb_2^2 \end{array}} && \gammacb_1^2 \rho^2 \Delta^2 + \gammacb_2^2 \\
\textrm{subject to} && \rho^2\Xcb \succeq (\lambdacb_1 \omega_1^2 + \lambdacb_2 \omega_1^2  \omega_2^2)\I_n, 
\\
&& \Xcb \succeq \I_n, \quad \Xcb \preceq \gammacb_1^2 \I_n \\
&&
\left[\begin{array}{c|c|c}
\begin{array}{cc} -\lambdacb_1 \I_n & 0 \\ &-\lambdacb_2 \I_n \end{array}  & {\bf 0}_{2n\times 2n} &  \begin{array}{c}\Xcb-\Wcb^\transp \\ \Wcb^\transp \end{array}\\  \hline
&   -\gammacb_2^2\I_{2n} & \begin{array}{c} \Q(\Xcb-\Wcb^\transp)\\ \R\Wcb^\transp \end{array} \\ \hline
 \star & & -\Xcb
\end{array} \right] \preceq 0,\label{condK}
\end{eqnarray}
\end{subequations}
%
%
then Algorithm~\ref{alg:contr} with $\K = \Wcb\Xcb^{-1}$ generates a sequence $\{\x_k\}_{k \in \naturals}$ that converges as,
\begin{equation}\label{ate2} 
\mathrm{AE} := \limsup_{k\to \infty} \E[\| \x_k - \hat{\x}^{\star}_k\|] \leq \frac{1}{1-\rho}\left(\gamma_1 \rho \Delta + \gamma_2 \right)\,. 
\end{equation}
Furthermore, for any fixed $\rho$, solving Problem~\eqref{condit.3} diminishes the asymptotical error $\mathrm{AE}$. \qed
\end{theorem}

Theorem~\ref{th.3} describes how to best tune the matrix $\K$ to minimize the asymptotic tracking error. In particular, doing a grid search on $\rho \in (0,1)$, one can identify the best $\K$ that minimizes the worst-case asymptotical error bound. Two remarks are in order. 

First the convex problem~\eqref{condit.3} grows linearly in the dimension of the problem $n$. However, if the matrices $\R$ and $\Q$ have some particular structure (i.e., diagonal, block diagonal), we can reduce the problem size considerably. In the case of $\Q = q \I_n$, $\R = r \I_n$, then choosing $\X = p\I_n, \W = w \I_n$, the problem becomes independent of the problem size $n$, as it happens in other examples~\cite{Lessard2022}. 

Second, the matrices $\R$ and $\Q$ are gateways through which one can include variable correlations and dependency, which is not present in a standard prediction-correction algorithm, nor in the scalar worst-case convergence tuning $\chi$.

\subsection{LPV filter design}

The presented filter design, captured by the conditions in Theorem~\ref{th.3}, could suffer from conservatism, since the matrices $\Q$ and $\R$ are selected as worst cases. For example, we have modeled $\Q$ at having its largest singular value at $\sqrt{2}(\xip \tau_{\mu})$. In practice, one could wish to model these matrices as parameter-varying, and directly dependent on how the data changes in time. 

We propose here an LPV filter design that accomplishes this task and reduce conservatism of the design process. We focus on a simplified setting to convey the basic ideas, and the reader is left to generalize the approach. 

Consider the change in time of the data $\nabla_t \y(t)$, and let $\theta \in [0,1]$ be a normalized parameter that capture this change. For example, we let $\theta = (\|\nabla_t \y(t)\|_{\infty})/(\max_{t} \|\nabla_t \y(t)\|_{\infty})$. We then consider $\Q$ as an affine parameter-varying matrix: $\Q(\theta) = \Q_0 + \theta \Q_1$ .  We consider here a static $\R$ for simplicity, thereby assuming that the change in the data only affects the prediction accuracy, which is reasonable to assume (and easy to lift if needed). 

In parallel, we will be looking for an affine parameter-varying Lyapunov matrix $\X(\theta) = \X_0 + \theta \X_1$, and a parameter-varying filter gain matrix $\K(\theta)$. The following theorem is then in place. 

\begin{theorem}\label{th.4}
Consider Algorithm~\ref{alg:contr} and its abstraction~\eqref{abstraction}, to find and track the filtered optimizer trajectory $\hat{\x}^{\star}(t)$ of the time-varying stochastic optimization problem $\min_{\x\in\reals^n} \E[f(\x;\y(t))]+g(\x)$. Assume that the optimizer trajectory varies in a bounded way as $\|\hat{\x}^{\star}(t_{k+1})- \hat{\x}^{\star}(t_k)\| \leq \Delta$, for all $k\in \naturals$ and $\Delta<\infty$. Let Assumptions~\ref{as.1}-\ref{as.2} hold as well.   

Introduce matrix function $\W(\theta): [0,1] \to \reals^{n\times n} = \W_0 + \theta \W_1$ and $\X(\theta): [0,1] \to \reals^{n\times n} = \X_0 + \theta \X_1, \X(\theta)\succ 0$, scalars $\lambda_1\geq0, \lambda_2\geq0$, in addition to supporting scalars $\gamma_1, \gamma_2$, and consider the tuning matrix function $\K(\theta)$. Assume that $\Q(\theta) = \Q_0 + \theta \Q_1$ and that the value of $\theta$ at subsequent time step is upper bounded as $|\theta_{s+1}- \theta_s|\leq \nu$. For any fixed scalar $\rho\in(0,1)$, by solving the problem
\begin{subequations}\label{condit.4}
\begin{eqnarray}
\hskip-1cm\minimize_{\scriptsize\begin{array}{c}\Xcb_0, \Xcb_1, \lambdacb_1\geq0, \lambdacb_2\geq0,\\ \Wcb_0, \Wcb_1, \gammacb_1^2, \gammacb_2^2 \end{array}} && \gammacb_1^2 \rho^2 \Delta^2 + \gammacb_2^2 \\
\textrm{subject to} && 
\rho^2[\Xcb_0 + \Xcb_1] \succeq (\lambdacb_1 \omega_1^2 + \lambdacb_2 \omega_1^2  \omega_2^2)\I_n, \label{eq.affine1} \\
&& [\Xcb_0 + \Xcb_1] \succeq \I_n, \quad \Xcb_0 \preceq \gammacb_1^2 \I_n \label{eq.affine2} \\
&& \Xcb_1 \preceq 0 
\label{eq.X1} \\
&&
\left.\begin{array}{l}\left[\begin{array}{c|c|c}
\begin{array}{cc} -\lambdacb_1 \I_n & 0 \\ &-\lambdacb_2 \I_n \end{array}  & {\bf 0}_{2n\times 2n} &  \begin{array}{c}{\Ycb}(\theta)-\Wcb(\theta)^\transp \\ \Wcb(\theta)^\transp \end{array}\\  \hline
&   -\gammacb_2^2\I_{2n} & \begin{array}{c} \Q(\theta)({\Ycb}(\theta)-\Wcb(\theta)^\transp)\\ \R\Wcb(\theta)^\transp \end{array} \\ \hline
 \star & & -{\Ycb}(\theta)
\end{array} \right] \preceq 0, \\ \\{\Ycb}(\theta) = \Xcb_0 - \nu \Xcb_1 + \theta \Xcb_1 \end{array}\right\} \quad \forall \theta \in [0,1] \label{eq.quadratic}
\end{eqnarray}
\end{subequations}
then Algorithm~\ref{alg:contr} with $\K(\theta) = \Wcb(\theta){\Ycb}(\theta)^{-1}$ generates a sequence $\{\x_k\}_{k \in \naturals}$ that converges as,
\begin{equation}\label{ate2} 
\mathrm{AE} := \limsup_{k\to \infty} \E[\| \x_k - \hat{\x}^{\star}_k\|] \leq \frac{1}{1-\rho}\left(\gamma_1 \rho \Delta + \gamma_2 \right)\,. 
\end{equation}
Furthermore, for any fixed $\rho$, solving Problem~\eqref{condit.3} diminishes the asymptotical error $\mathrm{AE}$. \qed
\end{theorem}

Theorem~\ref{th.4} describes parametric conditions for Algortihm~\ref{alg:contr} to converge to the optimizer trajectory in expectation and within an error ball. We remark here the extra constraint $\X_1 \preceq 0$, which requires some explanation. 

As one can see in the proof of Theorem~\ref{th.4}, the key step in proving convergence of the algorithm is to ensure that a Lyapunov function decreases at subsequent times. The Lyapunov function that we consider is $\mathcal{L}_k(\theta_k) = (\x_{k}-\hat{\x}^{\star}_k)^\top \X(\theta_k) (\x_{k}-\hat{\x}^{\star}_k)$, and therefore we have to deal with matrices $\X(\theta)$ at subsequent times. By imposing $\X_1 \preceq 0$, we can write however,
$$
\X(\theta_{s+1}) = \X(\theta_{s}) + (\theta_{s+1}-\theta_s) \X_1 \preceq \X(\theta_{s}) - |\theta_{s+1}-\theta_s| \X_1 = \Y(\theta_s),
$$
from which we can derive Theorem~\ref{th.4} and this renders the derivation of a solvable program easier. 

The constraint $\X_1 \preceq 0$ induces conservatism in the design, but in all our numerical simulations we found it to be redundant (meaning that the optimal $\X_1^{\star}$ was $\preceq 0$ with or without the constraint), and therefore not conservative in our application. We remark that a more correct approach would be to consider both extremes ($+\nu$ and $-\nu$) as done in~\cite{Massioni2020} and remove $\X_1 \preceq 0$, but this would \rev{require the introduction of four variables $\Y_{-}(\theta), \Y_{+}(\theta), \W_{-}(\theta), \W_{+}(\theta)$, and the determination of $\K(\theta) = \W_{-}(\theta)\Y_{-}(\theta)^{-1} \neq \W_{+}(\theta)\Y_{+}(\theta)^{-1} = \K(\theta)$ would be flawed; see the proof for more details}. 

Another possible approach is to remove $\X_1 \preceq 0$ and to consider only slowly changing parameters, for which $\nu \ll 1$. Since in our simulations $\nu \approx 0.4$, we have preferred to focus on the former approach. Finally, note that considering a $\X_1 \preceq 0$ is not totally unreasonable, since we can assume that increasing $\theta$, the convergence performance would be negatively affected. 

For fixed $\rho$, the problem to be solved is infinite dimensional (yet convex once $\theta$ is fixed). We recall that the constraints in~\eqref{eq.quadratic} are quadratic in $\theta$, due to the product with $\Q(\theta)$. 

A possible way to solve problem~\eqref{condit.4} is to discretize the domain with a uniform grid $\Theta :=\{0, \theta_1, \ldots, \theta_{q}, 1\}$ and impose~\eqref{eq.quadratic} for all the points of the grid. Another possible way is to introduce another variable $\eta = \theta^2 \in [0,1]$, and render affine~\eqref{eq.quadratic}. A more sophisticated way to proceed is to introduce a more conservative yet convex condition in $\theta$, hinging on the concept of matrix sum of squares. 

\begin{definition}
Let $P(\theta)$ be a symmetric matrix of polynomials of degree up to $2d \in \naturals$ in the variable $\theta \in \reals$. A matrix is sum of squares (MSOS) if there exists a finite number $l$ of symmetric matrices of polynomials $\Pi_i(\theta)$ such that 
$$
P(\theta) = \sum_{i=1}^l \Pi_i(\theta)^\transp \Pi_i(\theta).
$$
The decomposition implies that $P(\theta)\succeq 0$ for all $\theta$. The constraint ``$P(\theta)$ is MSOS'' is convex.
\end{definition}

We are now ready for the following result. 

\begin{theorem}\label{th:MSOS}
Consider the matrix multiplicator $\Lambdacb \in \reals^{5n \times 5n}, \Lambdacb \succeq 0$. In Theorem~\ref{th.4}, condition~\eqref{eq.quadratic} can be substituted with the more conservative, yet convex and finite dimensional condition:
\begin{subequations}
\begin{eqnarray}
&&\hspace*{-6mm}
-\left[\begin{array}{c|c|c}
\begin{array}{cc} -\lambdacb_1 \I_n & 0 \\ &-\lambdacb_2 \I_n \end{array}  & {\bf 0}_{2n\times 2n} &  \begin{array}{c}{\Ycb}(\theta)-\Wcb(\theta)^\transp \\ \Wcb(\theta)^\transp \end{array}\\  \hline
&   -\gammacb_2^2\I_{2n} & \begin{array}{c} \Q(\theta)({\Ycb}(\theta)-\Wcb(\theta)^\transp)\\ \R\Wcb(\theta)^\transp \end{array} \\ \hline
 \star & & -{\Ycb}(\theta)
\end{array} \right]  - \Lambdacb \theta (1-\theta) \quad \textrm{ is MSOS},  \label{condKT5} \\ && \Lambdacb \succeq 0, \quad {\Ycb}(\theta) = \Xcb_0 - \nu \Xcb_1 + \theta \Xcb_1
\end{eqnarray}
\end{subequations}
where $\Lambdacb$ is now a new decision variable in the optimization problem. \qed 
\end{theorem}

Theorems~\ref{th.4} and~\ref{th:MSOS} describe an LPV design strategy for our filter gain design. This is a particular choice due to the affine parameter-varying $\Q(\theta)$ and static $\R$. More complex choices can be made (relatively) straightforwardly following the same pattern of the presented theorems. We now look at some numerical simulations. 

\section{Numerical simulations}\label{sec:num}

We focus now on showcasing the performance of the proposed algorithms on a real dataset and problem stemming from ride-hailing. We obtain trips data from the New York City dataset\footnote{Open data from the NYC Taxi and Limousine Commission Data Hub.} for the yellow taxi cab in the month of November 2019, totalling over $6.8$ millions trips. We group the trips in $5$ minutes intervals and divide the trip requests among $n=5$ different ride-hailing companies (such as taxi cab, Uber, Lyft, etc.). The trip requests are divided randomly among the companies in a way that different companies do not have the same number of requests. 

In the modern context of mobility as a service with ride-hailing orchestration~\cite{rs2}, it makes sense for the city to provide software platforms to decide caps on the number of vehicles that each company can put on the streets depending on trading-off satisfying the demand and limiting traffic. A natural optimization problem that a city can formulate is
\begin{equation}\label{rs}
\hat{\x}^\star(t) = \arg\min_{\x\in[\underline{x}, \overline{x}]^n} \sum_{i=1}^n \Big({\E}\Big[\frac{1}{2}\|\x_i-c_i\y_i(t)\|^2\Big]+ \log(1 + \kappa \exp(\x_i)) + \frac{\varsigma}{2} \sum_{j=1}^n \|\x_i-\x_j\|^2 \Big),
\end{equation} 
where $\x_i$ for each $i$ represents the upper limit on vehicles on the roads for company $i$. The constraint $[\underline{x}, \overline{x}]^n$ represents box constraints on the number of allowed vehicles. The term $c_i>0$ multiplies the trip requests for company $i$ at time $t$, $\y_i(t)$ to be able to match most of the requests as possible. The logistic term with $\kappa>0$ is a regularization to make the cost non-quadratic but still convex and favor a smaller number of vehicles on the roads. Finally, the coupling term $\|\x_i-\x_j\|^2$ is set to have a similar regulation among different companies. 
For the sake of the simulations, we take $c_i =1, \kappa = 0.02, \varsigma = 0.1$, and we sample Problem~\eqref{rs} at every $5$ minutes. We simulate also the ground truth considering a smoothed version of the data. 

\subsection{Unconstrained case}

We start by considering the unconstrained case, where $\x \in \reals^5$, and we look at different noise regimes and algorithms.  

As for the algorithms, we consider our TV-EKF (Algorithm~\ref{alg:cap}), a standard prediction-correction~\cite{Bastianello2020}, the stochastic online algorithm of~\cite{Cutler2021}, and the stochastic prediction-correction version of~\cite{Maity2022}. For the latter, with their optimal choice of window size and weights for two-point evaluation, we can see that it is equivalent to the AGT algorithm of~\cite{Paper1}, i.e., exact prediction with Hessian inversion, but with a finite difference evaluation of $\nabla_{t\x}f$.  

For the four algorithms, we consider different choices of prediction steps $P$ and correction steps $C$. For the stochastic online and prediction-correction version of~\cite{Cutler2021, Maity2022}, prediction is either absent ($P=0$) or exact ($P \to \infty$), so $P$ is not a free parameter. 
We simulate three noise regimes:  
 
$\bullet$ The case of very good prediction $Q/R \approx 0$. For this, $\OJ$ generates as predictive signal a random signal around the ground truth with variance $10$. We use for the correction the true data stream. This case represents a realistic scenario of having a very good predictor (based on accurate historical data and, e.g., on a periodic Kernel method). This could be typical in ride-hailing systems. 

$\bullet$ The case of poor prediction $R/Q \approx 0$. For this, $\OJ$ generates as predictive signal a random signal around the ground truth with variance $200$. We use for the correction a convex combination of the true data stream and the ground truth (weighting the true data $0.05$). This case represents a potential scenario of situation in which the system transitions (e.g., a lock-down happens) and we have poor prediction. This is in general less typical, since prediction can be built online on current data, based, e.g., on extrapolation, but still interesting to analyze. 

$\bullet$ The case of $\OJ$ generated by the true data based on an extrapolation-based prediction with \rev{one, two, and three-points}~\cite{Bastianello2020}, and correction also based on true data. We label this case $Q\approx R$. The error between the true data and the ground truth has variance $\approx 50$, and the prediction $\approx 200$. We design $\Q_k$ and $\R_k$ accordingly, also taking into account that the data streams are correlated, and therefore $\Q_k, \R_k$ are full. 

Table~\ref{tab.1} displays the results obtained in these settings. As we can see, Algorithm~\ref{alg:cap} performs well in most cases by a significant margin, when prediction is accurate ($Q/R\approx 0$). When prediction is poor ($R/Q\approx 0$), then Algorithm~\ref{alg:cap} behaves close to a damped Newton's method, which significantly outperforms a standard prediction-correction algorithm, but it is in par with its stochastic version (since the latter still uses the very accurate data stream to built its exact prediction). 

Finally, for the setting of $Q\approx R$, then all the four methods perform very similarly, which is almost expected since taking prediction, or prediction and then correction incur the same ``error'', so any combination could achieve similar results. In this case, Algorithm~\ref{alg:cap} chooses a $\K_k \approx \I_n$. 

The results in Table~\ref{tab.1} support Algorithm~\ref{alg:cap} as an algorithm that can automatically tune prediction and correction; based on this tuning it can be significantly better than the competitors; and in the worst case it performs as state-of-the-art methods. We finally remark that the good prediction case is considered to be typical in this application scenario, and that the method from~\cite{Maity2022} could also be updated by designing an EKF for it, in the same way we did for the standard prediction-correction.

\begin{remark}
\rev{Different choices of number of prediction and correction steps $P$ and $C$ yield different errors, and a higher number of steps does not necessarily means a lower error. This is not in contrast with Equation~\eqref{ate.w} of Theorem~\ref{prop.2}, since the theoretical result is an expected worst-case error based on the $\tau$ and $\sigma$ worst-case bounds. In practice however, if, for example, prediction is very inaccurate, then converging to it may not be the best course of action for a particular realization, so a lower number of prediction steps (at the limit $P=0$ for~\cite{Cutler2021}) may incur in lower errors. } 
\end{remark}

\begin{remark}
\rev{Usually, the Extended Kalman filter may be sensitive to the choice of the initial condition. In these simulations, we have chosen $0$ as the initial condition, but we have also tried different random conditions in the range $[0,1000]$ without seeing any difference in performance. We expect this to be due to the fact that both prediction and correction are contractive operator that could converge separately; however, further investigations are left for future endeavors. } 
\end{remark}


\begin{table}\centering
\caption{Performance of the considered algorithms in an unconstrained setting. For each row, the first line represents the \rev{average error: $\|\x_k-\hat{\x}^\star_k\|$}, the second line the $25\%$ percentile, and the third line the $75\%$ percentile. In bold, the smallest error for the selected case and parameter choice. With $^*$ we indicate a $>10\%$ error reduction with respect to the closest competitor. \rev{Note that the algorithms of~\cite{Cutler2021} and \cite{Maity2022} are not affected by the choice of prediction.} }
\label{tab.1}
\scalebox{0.9}{
\begin{tabular}{l|cccc|cc|cc|cccc}
\toprule
Regime & \multicolumn{12}{c}{Algorithm} \\ \toprule
& \multicolumn{4}{c}{Extrap. P-C, \cite{Bastianello2020}, $(P,C)$} & \multicolumn{2}{c}{Stoch. C~\cite{Cutler2021}, $C$} & \multicolumn{2}{c}{Stoch. P-C~\cite{Maity2022}, $C$} & \multicolumn{4}{c}{TV-EKF: Algortihm~\ref{alg:cap}, $(P,C)$} \\ \midrule
& $(1,1)$ &$(5,1)$ &$(1,5)$ &$(5,5)$ & $1$ & $5$ & $1$ & $5$ & $(1,1)$ &$(5,1)$ &$(1,5)$ &$(5,5)$ \\ \toprule 
$Q/R \approx 0$ &80.0 & 47.2 & 134.0 & 118.9 & 130.6& 142.8& 160.1 & 151.1 & {\bf{70.1}$^*$} & {\bf{10.8}$^*$} & {\bf{61.2}$^*$} & {\bf{14.0}$^*$} \\ 
 &{\color{darkgray}\footnotesize 26.3 } & {\color{darkgray}\footnotesize 17.0 } & {\color{darkgray}\footnotesize 49.7 } & {\color{darkgray}\footnotesize 44.7 }& {\color{darkgray}\footnotesize 40.9 } & {\color{darkgray}\footnotesize 52.3 } & {\color{darkgray}\footnotesize 59.9 } & {\color{darkgray}\footnotesize 56.0 } & {\color{darkgray}\footnotesize 20.2 } & {\color{darkgray}\footnotesize 3.8 } & {\color{darkgray}\footnotesize 17.7 } & {\color{darkgray}\footnotesize 4.6 } \\ 
 &  {\color{darkgray}\footnotesize 111.3 } & {\color{darkgray}\footnotesize 64.7 } & {\color{darkgray}\footnotesize 183.9 } & {\color{darkgray}\footnotesize 163.5 }  & {\color{darkgray}\footnotesize 184.7  }& {\color{darkgray}\footnotesize 195.5 } & {\color{darkgray}\footnotesize 219.1 } & {\color{darkgray}\footnotesize 206.6 } & {\color{darkgray}\footnotesize 107.1 } & {\color{darkgray}\footnotesize 14.9 } & {\color{darkgray}\footnotesize 90.8 } & {\color{darkgray}\footnotesize 19.3 }\\ \midrule
$R/Q \approx 0$ &90.5 & 195.7 & 19.4 & 48.6& 72.2& 10.5 & 11.2 & 8.0 & {\bf{7.5}$^*$} & {\bf{7.6}$^*$} & {\bf{7.5}} & {\bf{7.5}} \\ 
 & {\color{darkgray}\footnotesize 66.3 } & {\color{darkgray}\footnotesize 148.7 } & {\color{darkgray}\footnotesize 12.7 } & {\color{darkgray}\footnotesize 32.2 } & {\color{darkgray}\footnotesize 21.4 }& {\color{darkgray}\footnotesize 3.6 } & {\color{darkgray}\footnotesize 4.1 } & {\color{darkgray}\footnotesize 2.9 } & {\color{darkgray}\footnotesize 2.8 } & {\color{darkgray}\footnotesize 2.9 } & {\color{darkgray}\footnotesize 2.8 } & {\color{darkgray}\footnotesize 2.8 } \\ 
 & {\color{darkgray}\footnotesize 109.3 } & {\color{darkgray}\footnotesize 236.4 } & {\color{darkgray}\footnotesize 23.8 } & {\color{darkgray}\footnotesize 60.6 } & {\color{darkgray}\footnotesize 110.9 }& {\color{darkgray}\footnotesize 14.5 } & {\color{darkgray}\footnotesize 14.3 } & {\color{darkgray}\footnotesize 10.6 } & {\color{darkgray}\footnotesize 10.4 } & {\color{darkgray}\footnotesize 10.4 } & {\color{darkgray}\footnotesize 10.5 } & {\color{darkgray}\footnotesize 10.5 } \\ \midrule
${\color{white} Q \approx R}\qquad \OJ^{(1)}$ &132.2 & { 140.2} & 143.3 & 144.1 &{\bf 130.6}  & {\bf 142.8}& 160.1 & 151.1 & {{139.3}} & {{140.9}} & {{149.1}} & {{149.1}} \\ 
 & {\color{darkgray}\footnotesize 45.8 } & {\color{darkgray}\footnotesize 49.7 } & {\color{darkgray}\footnotesize 52.5 } & {\color{darkgray}\footnotesize 53.2 } & 
{\color{darkgray}\footnotesize 40.9 } & {\color{darkgray}\footnotesize 52.3 }& {\color{darkgray}\footnotesize 59.9 } & {\color{darkgray}\footnotesize 56.0 } & {\color{darkgray}\footnotesize 49.8 } & {\color{darkgray}\footnotesize 51.3 } & {\color{darkgray}\footnotesize 56.2 } & {\color{darkgray}\footnotesize 56.1 } \\ 
 & {\color{darkgray}\footnotesize 181.6 } & {\color{darkgray}\footnotesize 192.1 } & {\color{darkgray}\footnotesize 196.4 } & {\color{darkgray}\footnotesize 197.3 } & 
{\color{darkgray}\footnotesize 184.7 } & {\color{darkgray}\footnotesize 195.5 }& {\color{darkgray}\footnotesize 219.1 } & {\color{darkgray}\footnotesize 206.6 } & {\color{darkgray}\footnotesize 190.8 } & {\color{darkgray}\footnotesize 192.8 } & {\color{darkgray}\footnotesize 204.5 } & {\color{darkgray}\footnotesize 204.2 } \\  
${\color{black} Q \approx R}\qquad \OJ^{(2)}$ &{ 135.0} & 162.4 & { 144.7} & { 149.5} &{\bf 130.6}  & {\bf 142.8} & 160.1 & 151.1 & {{141.7}} & {{152.4}} & {{149.5}} & {{150.2}} \\ 
 & {\color{darkgray}\footnotesize 47.8 } & {\color{darkgray}\footnotesize 60.9 } & {\color{darkgray}\footnotesize 53.6 } & {\color{darkgray}\footnotesize 55.5 } &
 {\color{darkgray}\footnotesize 40.9 } & {\color{darkgray}\footnotesize 52.3 }& {\color{darkgray}\footnotesize 59.9 } & {\color{darkgray}\footnotesize 56.0 } & {\color{darkgray}\footnotesize 52.0 } & {\color{darkgray}\footnotesize 56.4 } & {\color{darkgray}\footnotesize 56.4 } & {\color{darkgray}\footnotesize 56.4 } \\ 
 & {\color{darkgray}\footnotesize 185.7 } & {\color{darkgray}\footnotesize 225.9 } & {\color{darkgray}\footnotesize 198.0 } & {\color{darkgray}\footnotesize 204.3 } &
 {\color{darkgray}\footnotesize 184.7 } & {\color{darkgray}\footnotesize 195.5 }& {\color{darkgray}\footnotesize 219.1 } & {\color{darkgray}\footnotesize 206.6 } & {\color{darkgray}\footnotesize 193.4 } & {\color{darkgray}\footnotesize 210.6 } & {\color{darkgray}\footnotesize 205.0 } & {\color{darkgray}\footnotesize 206.0 } \\  
${\color{white} Q \approx R}\qquad \OJ^{(3)}$ &137.8 & 199.2 & 145.4 & 156.6 &{\bf 130.6}  & {\bf 142.8} & 160.1 & 151.1 & {{143.1}} & {{170.3}} & {{149.6}} & {{151.0}} \\ 
 & {\color{darkgray}\footnotesize 49.0 } & {\color{darkgray}\footnotesize 73.4 } & {\color{darkgray}\footnotesize 53.5 } & {\color{darkgray}\footnotesize 57.9 } & 
  {\color{darkgray}\footnotesize 40.9 } & {\color{darkgray}\footnotesize 52.3 }&{\color{darkgray}\footnotesize 59.9 } & {\color{darkgray}\footnotesize 56.0 } & {\color{darkgray}\footnotesize 52.9 } & {\color{darkgray}\footnotesize 62.8 } & {\color{darkgray}\footnotesize 56.3 } & {\color{darkgray}\footnotesize 56.6 } \\ 
 & {\color{darkgray}\footnotesize 190.3 } & {\color{darkgray}\footnotesize 279.0 } & {\color{darkgray}\footnotesize 198.6 } & {\color{darkgray}\footnotesize 217.6 } &
 {\color{darkgray}\footnotesize 184.7 } & {\color{darkgray}\footnotesize 195.5 }& {\color{darkgray}\footnotesize 219.1 } & {\color{darkgray}\footnotesize 206.6 } & {\color{darkgray}\footnotesize 196.1 } & {\color{darkgray}\footnotesize 237.5 } & {\color{darkgray}\footnotesize 205.5 } & {\color{darkgray}\footnotesize 206.7 } \\ \bottomrule
\end{tabular}}
\end{table}

\subsection{Constrained case}

We analyze now the constrained case, for which we set $\underline{x} = 100$ and $\overline{x} = 1000$. These constraints are not overly restrictive, but the point here is to see how our BMI-based method performs with respect to a standard prediction-correction method in different scenarios. Here, the stochastic variant~\cite{Maity2022} cannot be applied, but one could use the methods in~\cite{Simonetto2020} with exact prediction and finite-difference computations for $\nabla_{t\x}f$. We do not look into that, since typically these methods are more computationally demanding. 

The settings we investigate are similar to the unconstrained case, with the difference that we use our second algorithm TV-CONTRACT with a static $\K$ generated via Problem~\eqref{condit.3}, and uniform grid-search on $\rho$. We also consider two cases for $Q\approx R$ \rev{with a two-point extrapolation prediction}, the first \rev{$(a)$} has $\Q \approx 200\I_n$ and $\R \approx 50 \I_n$, the second \rev{$(b)$}: $\Q \approx 67\I_n$ and $\R \approx 50 \I_n$. In both cases, as before, $\Q, \R$ are full.    

In Table~\ref{tab.2}, we displays the results obtained in these settings. As we can see, Algorithm~\ref{alg:contr} has the best advantage when prediction is accurate and $\K$ can be chosen different from $\I_n$. In general, the performance of Algorithm~\ref{alg:contr} is comparable with the competition, unless there is a clear advantage to choose a different from $\I_n$ gain. In this latter case, the performance gain can be important. \rev{For the case $R/Q \approx 0$, i.e., when prediction is really poor, we select $\K \approx \I_n$ and we perform in part with~\cite{Bastianello2020}, but slightly worse than doing no prediction at all, as in~\cite{Cutler2021}.} In Table~\ref{tab.2}, we have also added the selected best $\K$, which is full whenever we use the $\approx$ sign, with diagonal elements close to the indicated values.

The results in Table~\ref{tab.2} support Algorithm~\ref{alg:contr} as an algorithm that can automatically tune prediction and correction; based on this tuning it can be better than the competitors; and in the worst case it performs in par with state-of-the-art methods. We finally remark that the good prediction case is considered to be typical in this application scenario.
 
\begin{figure}\centering
\resizebox{0.8\textwidth}{!}{\input{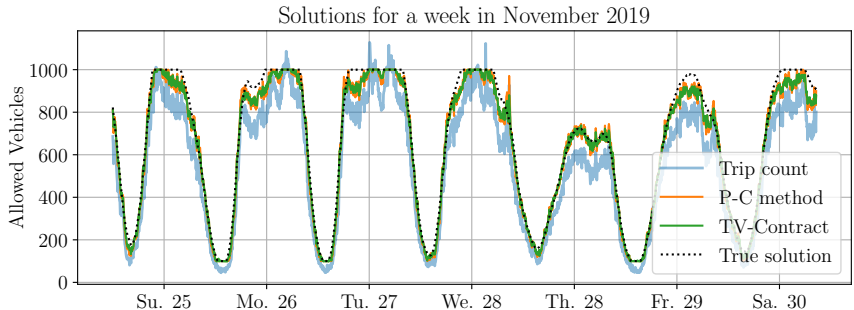}}
\caption{Display of selected trajectories in Thanksgiving week of 2019. The setting is $Q\approx R (b)$ and $P=C=5$, and the trajectories are for one of the five companies. }
\label{fig.3}
\end{figure}

For display purposes, Figure~\ref{fig.3} illustrates the different trajectories for one of the five companies during Thanksgiving week of the selected month of November 2019, for the case $Q\approx R (b)$ and $P=C=5$.

\begin{table}\centering
\caption{Performance of the considered algorithms in a constrained setting. For each row, the first line represents the \rev{average error: $\|\x_k-\hat{\x}^\star_k\|$}, the second line the $25\%$ percentile, and the third line the $75\%$ percentile. In bold, the smallest error for the selected case and parameter choice. With $^*$ we indicate a $>10\%$ error reduction with respect to the closest competitor. \rev{Note that the algorithms of~\cite{Cutler2021} is not affected by the choice of prediction.}}
\label{tab.2}
\scalebox{0.9}{
\begin{tabular}{c|cccc|cc|ccccc}
\toprule
Regime & \multicolumn{10}{c}{Algorithm} & $\K_{(5,5)}$ \\ \toprule
& \multicolumn{4}{c}{Extrapolation P-C, \cite{Bastianello2020}, $(P,C)$} & \multicolumn{2}{c}{Stoch. C~\cite{Cutler2021}, $C$} & \multicolumn{4}{c}{TV-CONTRACT: Algorithm~\ref{alg:contr}, $(P,C)$}& \\
& $(1,1)$ &$(5,1)$ &$(1,5)$ &$(5,5)$ & $1$& $5$ & $(1,1)$ &$(5,1)$ &$(1,5)$ &$(5,5)$& \\ \toprule 
$Q/R \approx 0$ &{\bf 68.7} & 58.2 & 84.3 & 79.7 & 88.3 & 87.0 & 71.1 & {\bf 51.3$^*$} & {\bf 75.4$^*$} & {\bf 54.6$^*$} & \\
& {\color{darkgray}\footnotesize  25.3 }& {\color{darkgray}\footnotesize  20.3 }& {\color{darkgray}\footnotesize  33.5 }& {\color{darkgray}\footnotesize  32.6 }& {\color{darkgray}\footnotesize  28.5 }& {\color{darkgray}\footnotesize  33.8 }& {\color{darkgray}\footnotesize  32.8 }& {\color{darkgray}\footnotesize  9.2 }& {\color{darkgray}\footnotesize  30.5 }& {\color{darkgray}\footnotesize  15.3 } & $\K \approx 0.24\I_n $ \\ 
& {\color{darkgray}\footnotesize  102.4 }& {\color{darkgray}\footnotesize  93.1 }& {\color{darkgray}\footnotesize  122.1 }& {\color{darkgray}\footnotesize  117.1 }& {\color{darkgray}\footnotesize  119.1 }& {\color{darkgray}\footnotesize  125.1 }& {\color{darkgray}\footnotesize  103.1 }& {\color{darkgray}\footnotesize  93.8 }& {\color{darkgray}\footnotesize  108.5 }& {\color{darkgray}\footnotesize  92.7 } &  \\ \midrule 
$R/Q \approx 0$ &88.4 & 142.6 & 54.4 & 66.5 & \bf 83.8 & \bf 51.3 & {88.4} & {142.6} & { 54.4} & { 66.5} & \\
& {\color{darkgray}\footnotesize  57.4 }& {\color{darkgray}\footnotesize  79.7 }& {\color{darkgray}\footnotesize  16.0 }& {\color{darkgray}\footnotesize  33.3 }& {\color{darkgray}\footnotesize  43.9 }& {\color{darkgray}\footnotesize  9.4 }& {\color{darkgray}\footnotesize  57.4 }& {\color{darkgray}\footnotesize  79.7 }& {\color{darkgray}\footnotesize  16.0 }& {\color{darkgray}\footnotesize  33.3 } & $\K = \I_n$ \\ 
& {\color{darkgray}\footnotesize  120.0 }& {\color{darkgray}\footnotesize  199.4 }& {\color{darkgray}\footnotesize  95.8 }& {\color{darkgray}\footnotesize  98.8 }& {\color{darkgray}\footnotesize  113.6 }& {\color{darkgray}\footnotesize  94.4 }& {\color{darkgray}\footnotesize  120.0 }& {\color{darkgray}\footnotesize  199.4 }& {\color{darkgray}\footnotesize  95.8 }& {\color{darkgray}\footnotesize  98.8 } &  \\ \midrule 
$Q \approx R $ & 85.4 & \bf 93.3 & 87.7 & 89.2 &  88.3 & \bf 87.0 & {\bf 85.3} & { 94.2} & { 87.7} & { 89.2} & \\
$(a)$ & {\color{darkgray}\footnotesize  32.5 }& {\color{darkgray}\footnotesize  34.4 }& {\color{darkgray}\footnotesize  34.2 }& {\color{darkgray}\footnotesize  34.9 }& {\color{darkgray}\footnotesize  28.5 }& {\color{darkgray}\footnotesize  33.8 }& {\color{darkgray}\footnotesize  32.3 }& {\color{darkgray}\footnotesize  34.4 }& {\color{darkgray}\footnotesize  34.2 }& {\color{darkgray}\footnotesize  34.9 } & $\K \approx \I_n$ \\ 
& {\color{darkgray}\footnotesize  118.8 }& {\color{darkgray}\footnotesize  135.6 }& {\color{darkgray}\footnotesize  126.3 }& {\color{darkgray}\footnotesize  128.8 }& {\color{darkgray}\footnotesize  119.1 }& {\color{darkgray}\footnotesize  125.1 }& {\color{darkgray}\footnotesize  118.7 }& {\color{darkgray}\footnotesize  137.7 }& {\color{darkgray}\footnotesize  126.3 }& {\color{darkgray}\footnotesize  128.8 } &  \\ \midrule 
$Q \approx R $ &\bf 68.3 & 59.0 & \bf 84.1 & 79.3 & 88.3 & 87.0 & {\bf 68.3} & {\bf 56.8} & {\bf 84.1} & {\bf 68.8$^*$} & \\
$(b)$& {\color{darkgray}\footnotesize  25.0 }& {\color{darkgray}\footnotesize  22.8 }& {\color{darkgray}\footnotesize  33.5 }& {\color{darkgray}\footnotesize  32.4 }& {\color{darkgray}\footnotesize  28.5 }& {\color{darkgray}\footnotesize  33.8 }& {\color{darkgray}\footnotesize  27.4 }& {\color{darkgray}\footnotesize  21.5 }& {\color{darkgray}\footnotesize  33.5 }& {\color{darkgray}\footnotesize  28.3 } & $\K \approx 0.86 \I_n $ \\ 
& {\color{darkgray}\footnotesize  102.1 }& {\color{darkgray}\footnotesize  93.1 }& {\color{darkgray}\footnotesize  122.0 }& {\color{darkgray}\footnotesize  116.6 }& {\color{darkgray}\footnotesize  119.1 }& {\color{darkgray}\footnotesize  125.1 }& {\color{darkgray}\footnotesize  102.1 }& {\color{darkgray}\footnotesize  92.6 }& {\color{darkgray}\footnotesize  122.0 }& {\color{darkgray}\footnotesize  103.3 } &  \\ \bottomrule
\end{tabular}}
\end{table}

\subsection{Variable $K$ case}

We finish the simulation assessment showing the tracking results obtained solving Problem~\eqref{condit.4} for an affine parametric-varying $\Q$. In particular, we let $\Q_0 = \Q/5$ and $\Q_1 = 4 \Q/5$, where $\Q$ is numerically defined as before, and run Algorithm~\ref{alg:contr} on all the four cases that we have looked at in Table~\ref{tab.2}. We solve Problem~\eqref{condit.4} by uniform gridding with $4$ points. As mentioned, in our case $\nu \approx 0.4$. 

In Table~\ref{tab.3}, we report the results for the $Q \approx R$ cases, since we do not observe any substantial difference for the other cases of Table~\ref{tab.2}. We indicate in bold if we have a gain w.r.t. a static gain, and with a dagger, if we have also a gain w.r.t. the state of the art. We also report how the maximal element of the diagonal of $\K$ changes in time in a selected week. 

As we can see, the results are very similar to the static results, but the gains can be important in some cases. As we see, the filter gain does not change over a wide range. However, it appears that even these small changes are enough to reduce the asymptotical error in selected scenarios, and behaving in par with the static approach in the others. As one can infer, parametric-varying gain design does depend on the modelling choices for $\Q(\theta)$ and $\X(\theta)$, and one could expect possibly more performant results in the case of more complex dependencies on $\theta$. We leave this analysis for future endeavors.

\begin{table}\centering
\caption{Performance of the considered algorithm in a constrained setting. For each row, the first line represents the \rev{average error: $\|\x_k-\hat{\x}^\star_k\|$}, the second line the $25\%$ percentile, and the third line the $75\%$ percentile. We indicate in bold if we have a gain w.r.t. a static gain of Table~\ref{tab.2}, and with a dagger, if we have also a gain w.r.t. the state of the art of Table~\ref{tab.2}. Finally, with $^*$ we indicate a $>10\%$ error reduction with respect to the closest competitor. }
\label{tab.3}
\scalebox{0.9}{
\begin{tabular}{cccccc}
\toprule
Regime & \multicolumn{5}{c}{Algorithm}  \\ \toprule
&   \multicolumn{5}{c}{TV-CONTRACT-LPV: Algorithm~\ref{alg:contr}, $(P,C)$} \\
&  $(1,1)$ &$(5,1)$ &$(1,5)$ &$(5,5)$ &\\ \toprule 
$Q \approx R $ &  { 85.5} & {\bf 93.9} & {\bf 87.3} & { 89.8} &  \multirow{3}{*}{\resizebox{0.425\textwidth}{!}{\input{Figures/variable-gain-2.pgf}}}\\
$(a)$& {\color{darkgray}\footnotesize 32.0 }& {\color{darkgray}\footnotesize 34.3 }& {\color{darkgray}\footnotesize 33.7 }& {\color{darkgray}\footnotesize 35.1 } & \\ 
&{\color{darkgray}\footnotesize 118.6 }& {\color{darkgray}\footnotesize 136.7 }& {\color{darkgray}\footnotesize 125.4 }& {\color{darkgray}\footnotesize 129.2 }  & \\ 
&&&&&\\
&&&&&\\
&&&&&\\
&&&&&\\ \midrule 
$Q \approx R $ & { 70.9} & { 57.0} & {\bf 74.4}$^{*\dag}$ & {\bf 65.4}$^\dag$  & \multirow{3}{*}{\resizebox{0.425\textwidth}{!}{\input{Figures/variable-gain-3.pgf}}}\\
$(b)$& {\color{darkgray}\footnotesize 33.0 }& {\color{darkgray}\footnotesize 22.6 }& {\color{darkgray}\footnotesize 29.2 }& {\color{darkgray}\footnotesize 26.5 }&\\ 
&  {\color{darkgray}\footnotesize 102.4 }& {\color{darkgray}\footnotesize 93.6 }& {\color{darkgray}\footnotesize 108.8 }& {\color{darkgray}\footnotesize 99.1 } & \\
&&&&&\\
&&&&&\\
&&&&&\\
&&&&&\\
\bottomrule 
\end{tabular}}
\end{table}

\section{Conclusions}
We have discussed several methods to generalize time-varying optimization algorithms to the case of noisy data streams. The methods are rooted in the intuition that prediction and correction can be seen as a nonlinear dynamical system and a nonlinear measurement equation, respectively. This leads to extended Kalman filter formulations as well as contractive filters based on bilinear matrix inequalities (BMI's). Numerical results are promising, even when using possibly conservative BMI conditions. 

\appendix

\section{Proofs}

\subsection{An additional example}
\begin{example}[Deterministic example]\label{ex:1}
Consider a deterministic method with an extrapolation predictor \rev{with two points}~\cite{Bastianello2020}, meaning: $\OJ_{k+1}(\x) = 2 \nabla_{\x} f(\x; \y_{k}) - \nabla_{\x} f(\x; \y_{k-1})$, where now $\y(t)$ is deterministic. 
Assume $f(\x;\y(t))$ is strongly convex and smooth, uniformly in $\y$, assume that the Hessian of $f(\x;\y(t))$ does not depend on $\y$, and assume the following bounds on data and mixed derivatives:
\begin{equation}\label{bounds}
\max\{\left\|\nabla_{t} \y(t)\right\|, \left\|\nabla_{tt} \y(t)\right\|, \left\|\nabla_{\y \x}  f(\x; \y(t))  \right\|,  \left\|\nabla_{\y \y \x}  f(\x; \y(t))  \right\|\} \leq C, \qquad \forall \x, t.
\end{equation}
Then we have that 
$$\|\OJ_{k+1}(\x^{\star}_{k+1}) - \nabla_{\x} f(\x^{\star}_{k+1}; \y_{k+1})\| \leq (C^2 + C^3) h^2.$$ \qed
\end{example}

\begin{proof} From~\cite[Proof of Lemma~3]{Bastianello2020}, we know that there exists a $\tau \in [t_{k-1}, t_{k+1}]$ such that
\begin{equation}
\left\|\OJ_{k+1}(\x^{\star}_{k+1}) - \nabla_{\x} f(\x^{\star}_{k+1}; \y_{k+1})  \right\| \leq \left\|\nabla_{tt} \nabla_{\x} f(\x^{\star}_{k+1}; \y(\tau)) \right\|\, h^2. 
\end{equation}
Let the $i$-th component of $\nabla_{\x} f(\x; \y(t))$ be $D_i(\x; \y(t))$. 

By using the higher-derivatives Fa\`a di Bruno's chain rule:

\begin{multline}\label{result.bruno}
\left[\nabla_{tt} \nabla_{\x} f(\x; \y(t))\right]_{i} = \frac{\partial^2}{\partial t^2} D_i(\x; \y(t))  = \sum_j \left(\frac{\partial}{\partial y_j} D_i(\x; \y(t))\, \frac{\partial^2 y_j(t)}{\partial t^2}\right) +\\ +\sum_{j, \ell} \left( \frac{\partial^2}{\partial y_j \partial y_\ell}D_i(\x; \y(t))\, \frac{\partial y_j(t)}{\partial t} \, \frac{\partial y_\ell(t)}{\partial t}\right),
\end{multline}
from which the thesis follows. 
\end{proof}

\subsection{Derivations for Example~\ref{ex:2}}\label{ap.2}

We choose to write $\x^{\star}= \x^{\star}_{k+1}$ as a short-hand notation, in this proof only. 

By linearity of $\nabla_{\x}f(\x^{\star}; \y(t))$ with respect to the parameter $\y(t)$, and the linearity of the expectation, \rev{we can write the following.

{\bf{Predictor $\OJ^{(1)}_{k+1}$.} }

We have,
\begin{eqnarray*}
&&\E_{\w \in \mathcal{Y}_{k}}[\|\nabla_{\x} f(\x^{\star}; \w) - \E_{\y \in \mathcal{Y}_{k+1}}[\nabla_{\x} f(\x^{\star}; \y)]\|]=\E_{\w \in \mathcal{Y}_{k}}[\|\nabla_{\x} f(\x^{\star}; \w  - \E_{\y \in \mathcal{Y}_{k+1}}[\y])\|] \\ &&= \E_{\w \in \mathcal{Y}_{k}}[\|\nabla_{\x} f(\x^{\star}; \w - \bar{\y}_{k+1})\|] = \E_{\e^1 \in \mathcal{N}(0, \Sigmab_{k})}[\|\nabla_{\x} f(\x^{\star}; \bar{\y}_{k} - \bar{\y}_{k+1}) + \nabla_{\x} f(\x^{\star}; \e^1 )\|].
\end{eqnarray*}
We now use the Triangle inequality, the fact that $\|\nabla_{\x} f(\x^{\star}; \bar{\y}_{k} - \bar{\y}_{k+1})\|\leq C_0 C h$ by linearity and the assumptions, and the mean value theorem for the nominal trajectory, to upper bound the last inequality as 
\begin{eqnarray*}
&&\E_{\e^1 \in \mathcal{N}(0, \Sigmab_{k}))}[\|\nabla_{\x} f(\x^{\star}; \bar{\y}_{k}- \bar{\y}_{k+1})\| + \|\nabla_{\x} f(\x^{\star};\e^1)\|]\\
&&=C_0 C h + \E_{\e^1 \in \mathcal{N}(0, \Sigmab_{k})}[\|\nabla_{\x} f(\x^{\star}; \e^1)\|]\leq C_0 C h^2 + C_0\, \E_{\e^1 \in \mathcal{N}(0, \Sigmab_{k})}[\|\e^1\|]  \leq C_0 C h + C_0 \Sigma.
\end{eqnarray*}

{\bf{Predictor $\OJ^{(2)}_{k+1}$} }

We have,}
\begin{eqnarray*}
&&\E_{\w \in \mathcal{Y}_{k}, \z \in \mathcal{Y}_{k-1}}[\|\nabla_{\x} f(\x^{\star}; 2\w - \z) - \E_{\y \in \mathcal{Y}_{k+1}}[\nabla_{\x} f(\x^{\star}; \y)]\|]\\&&= \E_{\w \in \mathcal{Y}_{k}, \z \in \mathcal{Y}_{k-1}}[\|\nabla_{\x} f(\x^{\star}; 2\w - \z - \E_{\y \in \mathcal{Y}_{k+1}}[\y])\|] \\ &&= \E_{\w \in \mathcal{Y}_{k}, \z \in \mathcal{Y}_{k-1}}[\|\nabla_{\x} f(\x^{\star}; 2\w - \z - \bar{\y}_{k+1})\|] \\&&= \E_{\e^1 \in \mathcal{N}(0, \Sigmab_{k}), \e^2 \in \mathcal{N}(0, \Sigmab_{k-1})}[\|\nabla_{\x} f(\x^{\star}; 2\bar{\y}_{k} - \bar{\y}_{k-1} - \bar{\y}_{k+1}) + \nabla_{\x} f(\x^{\star}; 2\e^1 - \e^2)\|].
\end{eqnarray*}
We now use the Triangle inequality, the result~\eqref{result.bruno}, and the mean value theorem for the nominal trajectory, to upper bound the last inequality as 
\begin{eqnarray*}
&&\E_{\e^1 \in \mathcal{N}(0, \Sigmab_{k}), \e^2 \in \mathcal{N}(0, \Sigmab_{k-1})}[\|\nabla_{\x} f(\x^{\star}; 2\bar{\y}_{k} - \bar{\y}_{k-1} - \bar{\y}_{k+1})\| + \|\nabla_{\x} f(\x^{\star}; 2\e^1 - \e^2)\|]\\
&&=\|\nabla_{\x} f(\x^{\star}; 2\bar{\y}_{k} - \bar{\y}_{k-1} - \bar{\y}_{k+1})\| + \E_{\e^1 \in \mathcal{N}(0, \Sigmab_{k}), \e^2 \in \mathcal{N}(0, \Sigmab_{k-1})}[\|\nabla_{\x} f(\x^{\star}; 2\e^1 - \e^2)\|]\\
&&\leq C_0 C h^2 + C_0\, \E_{\e^1 \in \mathcal{N}(0, \Sigmab_{k}), \e^2 \in \mathcal{N}(0, \Sigmab_{k-1})}[\|2\e^1 - \e^2\|]  \leq C_0 C h^2 + 3 C_0 \Sigma.
\end{eqnarray*} 

\rev{
{\bf{Predictor $\OJ^{(3)}_{k+1}$} }

\rev{Similarly as before, we obtain,
\begin{eqnarray*}
&&\E_{\w \in \mathcal{Y}_{k}, \z \in \mathcal{Y}_{k-1}, \q \in \mathcal{Y}_{k-2}}[\|\nabla_{\x} f(\x^{\star}; 3\w - 3\z + \q) - \E_{\y \in \mathcal{Y}_{k+1}}[\nabla_{\x} f(\x^{\star}; \y)]\|] \\&&= \E_{\e^1 \in \mathcal{N}(0, \Sigmab_{k}), \e^2 \in \mathcal{N}(0, \Sigmab_{k-1}), \e^3\in \mathcal{N}(0, \Sigmab_{k-2})}[\|\nabla_{\x} f(\x^{\star}; 3\bar{\y}_{k} - 3\bar{\y}_{k-1} + \bar{\y}_{k-2} - \bar{\y}_{k+1}) +\\&&\qquad + \nabla_{\x} f(\x^{\star}; 3\e^1 - 3\e^2 + \e^3)\|].
\end{eqnarray*}
We now use the Triangle inequality, again~\cite[Proof of Lemma~3]{Bastianello2020}}, and the mean value theorem for the nominal trajectory, to upper bound the last inequality as 
\begin{eqnarray*}
&&\E_{\e^1 \in \mathcal{N}(0, \Sigmab_{k}), \e^2 \in \mathcal{N}(0, \Sigmab_{k-1}), \e^3\in \mathcal{N}(0, \Sigmab_{k-2})}[\|\nabla_{\x} f(\x^{\star}; 3\bar{\y}_{k} - 3\bar{\y}_{k-1} + \bar{\y}_{k-2}- \bar{\y}_{k+1})\| + \\&&\qquad+\|\nabla_{\x} f(\x^{\star}; 3\e^1 - 3\e^2+\e^3)\|]
\leq C_0 C h^3 + 7 C_0 \Sigma,
\end{eqnarray*} 
}
from which the first claim is proven.

For the second,
\begin{eqnarray*}
&&\E_{\y\in\mathcal{Y}_{k+1}} [\|\nabla_{\x}f(\x; \y) - \E_{\y \in \mathcal{Y}_{k+1}}[\nabla_{\x} f(\x; \y)] \|] =\E_{\y\in\mathcal{Y}_{k+1}} [\|\nabla_{\x}f(\x; \y) - \nabla_{\x} f(\x; \bar{\y}_{k+1})] \|] = \\ &&= \E_{\e\in\mathcal{N}(0,\Sigmab_{k+1})} [\|\nabla_{\x}f(\x; \e)] \|] \leq C_0 \, \E_{\e\in\mathcal{N}(0,\Sigmab_{k+1})} [\|\e \|] \leq C_0 \Sigma, 
\end{eqnarray*} 
as claimed.
\qed

\subsection{Proof of Proposition~\ref{pr.kalman}}\label{ap:pr.kalman}

Consider $C=1$ in Algorithm~\ref{alg:cap}, as well as a negligible $\R_k$. Then, we can simplify the Kalman gain as:
\begin{eqnarray*}
\K_k &=& \P_{k|k-1}\H_{k} [ \H_{k} \P_{{k|k-1}} \H_{k}^\transp]^{-1} =  \H_{k}^{-1}. 
\end{eqnarray*}
Therefore, the state update reads
\begin{eqnarray*}
\x_{k} &=& \x_{{k|k-1}} + \K_k (\Psib(\x_{k|k-1}, {\y}_{k})) = \x_{k|k-1} + \H_{k}^{-1}(-\x_{k|k-1} + \x_{k|k-1} - \beta \nabla_{\x} f(\x_{k|k-1}; {\y}_{k})) \\ &=& \x_{k|k-1} - \beta[\nabla_{\x\x}f(\x; {\y}_{k})]^{-1} \nabla_{\x} f(\x_{k|k-1}; {\y}_{k}),
\end{eqnarray*}
from which the thesis is proven. \qed

\subsection{Supporting results for Theorem~\ref{prop.2}}

\begin{lemma}\label{lemma.1}
Let Assumptions~\ref{as.1}-\ref{as.2} hold. Choose $\alpha<2\mu/L^2$. Let $\x_{k+1}^{\mathrm{f}}$ be the fixed point of the prediction ``pseudo''-dynamical model: $\x_{k+1}^{\mathrm{f}} = \Phib_{k,g}(\x_{k+1}^{\mathrm{f}})$. Then the distance between $\x_{k+1}^{\mathrm{f}}$ and the optimizer trajectory is bounded in expectation as,
\begin{equation*}
\E[\|\x_{k+1}^{\mathrm{f}} - \hat{\x}^{\star}_{k+1} \|] \leq  \frac{1}{\mu}\E [\|\OJ_{k+1}(\hat{\x}^{\star}_{k+1}) - \E_{\y \in \mathcal{Y}_{k+1}}[\nabla_{\x} f(\hat{\x}^{\star}_{k+1}; \y)]\|] \leq \frac{\tau}{\mu} =: \tau_{\mu}.\end{equation*}

Furthermore, let the setting of Example~\ref{ex:2} \rev{predictor $(2)$} hold. Then,
$$
\E[\|\x_{k+1}^{\mathrm{f}} - \hat{\x}^{\star}_{k+1} \|] \leq C_0 Ch^2/\mu + 3C_0\Sigma/\mu.
$$   \qed
\end{lemma}

\begin{proof} Choosing $\alpha<2\mu/L^2$ and under Assumption~\ref{as.1}, we know that the prediction is a contractive operator and its fixed point exists and it is unique. By implicit function theorems, see for instance~\cite[Theorem~2F.9]{Dontchev2009} and \cite[Theorem~3 and Lemma~2]{Bastianello2020}, being careful to $\OJ$ being strongly monotone and not generally the gradient of a strongly convex function, then,
\begin{equation}\label{diamond}
\|\x_{k+1}^{\mathrm{f}} - \hat{\x}^{\star}_{k+1}\| \leq \frac{1}{\mu} \underbrace{\|\OJ_{k+1}(\hat{\x}^{\star}_{k+1}) - \E_{\y \in \mathcal{Y}_{k+1}}[\nabla_{\x} f(\hat{\x}^{\star}_{k+1}; \y)]\|}_{(\diamond)}.
\end{equation}
Passing in expectations, and by using Assumption~\ref{as.2}, the first thesis follows. 

As for the second statement, it follows from the derivations of Example~\ref{ex:2} \rev{predictor $(2)$}. 
\end{proof}

\begin{lemma}\label{lemma.2}
Let Assumptions~\ref{as.1}-\ref{as.2} hold. Choose $\alpha<2\mu/L^2, \beta<2/L$. 
Consider the prediction update $\x_{k+1|k} = \Phib_{k,g}(\x_k)$ with $P$ prediction steps, and the correction update $\x_k' = \Psib_g'(\x_{k+1|k}, \y_{k+1})$ with $C$ correction steps. Let the contraction factors $\rhop, \rhoc$ be defined as in~\eqref{contra.factors}. Then, the following error bounds are in place. 
\begin{eqnarray}
\E[\|\x_{k+1|k} - \hat{\x}^{\star}_{k+1}\|] &\leq& \rhop^P \,\E[\|\x_{k} - \x_{k+1}^{\mathrm{f}}\|] + \tau_{\mu}\\ 
\E[\|\x_{k}' - \hat{\x}^{\star}_{k+1}\|] &\leq& \rhoc^C \,\E[\|\x_{k+1|k} - \hat{\x}^{\star}_{k+1}\|] + \sigmac,\qquad
\end{eqnarray}
where $\sigmac = \frac{\beta \,\sigma}{1-\rhoc}$.

Furthermore, under the setting of Example~\ref{ex:2}, \rev{predictor $(2)$}, $\tau_{\mu} = (C_0 Ch^2 + 3C_0\Sigma)/\mu$ and $\sigma = C_0 \Sigma$. \qed
\end{lemma}

\begin{proof}
Choosing $\alpha<2\mu/L^2, \beta<2/L$ and under Assumption~\ref{as.1}, we know that the prediction and correction are contractive operators and their fixed points are unique. 

For the prediction part, by using Equation~\eqref{diamond}, we obtain,
\begin{multline}
\|\x_{k+1|k} - \hat{\x}^{\star}_{k+1}\| =  \|\x_{k+1|k} \pm \x^\mathrm{f}_{k+1} - \hat{\x}^{\star}_{k+1}\| \leq \|[\prox_{\alpha g}(\OI - \alpha \OJ_{k+1}(\bullet))]^{\circ P}\x_{k} - \x^\mathrm{f}_{k+1}\| + \frac{1}{\mu} (\diamond) \leq \\ \leq \rhop^P \,\|\x_{k} - \x^\mathrm{f}_{k+1}\| + \frac{1}{\mu} (\diamond), \qquad
\end{multline} 
and passing in expectation with Assumption~\ref{as.2} the claim is proven. 

For the second claim, we can write
\begin{eqnarray}
\|\x_{k}' - \hat{\x}^{\star}_{k+1}\| \leq \|[\prox_{\beta g}(\OI - \beta \nabla_{\x}f(\bullet; \y_{k+1}) \pm \beta \E_{\y \in \mathcal{Y}_{k+1}}[\nabla_{\x} f(\bullet; \y)])]^{\circ C}\x_{k+1|k} - \hat{\x}^{\star}_{k+1}\|.\quad
\end{eqnarray}
Call $\epsilon_{k+1}(\x) := \nabla_{\x}f(\x; \y_{k+1}) - \E_{\y \in \mathcal{Y}_{k+1}}[\nabla_{\x} f(\x; \y)]$. Then each proximal gradient step will incur in an additive $\|\epsilon_{k+1}(\x)\|$ error, where $\x$ will be different at each step:
\begin{equation}
\|\x_{k}' - \hat{\x}^{\star}_{k+1}\| \leq \rhoc \|{\x}^{C-1}_{k+1|k} -\hat{\x}^{\star}_{k+1}\| + \beta\|\epsilon_{k+1}({\x}^{C-1}_{k+1|k})\|
\leq \rhoc^C \,\|\x_{k+1|k} - \hat{\x}^{\star}_{k+1}\| +  \underbrace{\beta \sum_{c=1}^C \rhoc^{C-c} \| \epsilon_{k+1}(\x^{C-c}_{k+1|k})\|}_{(\diamond\diamond)}. 
\end{equation}
Passing in expectation, with Assumption~\ref{as.2} and the sum of geometric series, the second claim is also proven. \end{proof}

\subsection{Proof of Theorem~\ref{prop.2}}\label{ap.prop.2}

The proof follows the one of \cite[Proposition 1]{Bastianello2020}, combining Lemma~\ref{lemma.1} and Lemma~\ref{lemma.2}. 
Start by considering $\chi = 1$, so a classical prediction-correction method. We can use \cite[Proposition 1]{Bastianello2020}, with $\E[\tau_k] = \tau_{\mu}$, and the correction with an additional error term to say,
\begin{equation}\label{dummy.01}
	\|{\x}_{k+1} - \hat{\x}^\star_{k+1}\| \leq  \zetac \Big(\zetap\|{{\x}_{k} - \hat{\x}^\star_{k}}\| + \zetap \Delta + \xip \tau_k \Big) + (\diamond\diamond) =: E_1.
\end{equation} 

Looking at prediction only, $\chi = 0$, we obtain instead,
\begin{equation}\label{dummy.02}
	\|{\x}_{k+1} - \hat{\x}^\star_{k+1}\| \leq  \Big(\zetap\|{{\x}_{k} - \hat{\x}^\star_{k}}\|   + \zetap \Delta + \xip \tau_k \Big)=: E_2.
\end{equation} 
 
For a generic $\chi \in [0,1]$, we can combine the errors as
\begin{equation}
	\|{\x}_{k+1} - \hat{\x}^\star_{k+1}\| \leq  (1-\chi)E_2 + \chi E_1.
\end{equation} 
Then, we can recursively compute the error via geometric series summation. By passing through expectations, the claim follows. 

For Example~\ref{ex:2}, with $\nabla_{\y\x}f$ bounded, by implicit function theorems~\cite{Simonetto2020}, we have that $\Delta = C_0 h/\mu$, from which the thesis.
\qed

\subsection{Proof of Proposition~\ref{lemma.3}}\label{ap.lemma.3}

For Problem~\eqref{p.w}, we look at the minimum of the curve,
\begin{equation}
\min_{\chi \in [0,1]}\,\frac{a \chi + b}{c \chi +d} =: F(\chi).
\end{equation}
For our problem $c \chi +d = 1-\zetap + \chi (\zetap-\zetap\zetac)>0$, $b = \zetap\Delta + \xip\tau_\mu>0$, while $a = (\zetac-1)(\zetap\Delta + \xip\tau_\mu) + \zetac\sigmac$ can be positive, negative, or zero. Since function $F(\chi)$ is a linear-fractional function in one dimension, for $\chi\geq0$, function $F(\chi)$ is monotone. In particular, for $a/c< (>) b/d$ the function is decreasing (increasing), leading to the optimal choices of $\chi^\star = 1 (0)$. The condition means,
\begin{equation*}
\frac{(\zetac-1)(\zetap\Delta + \xip\tau_\mu) + \zetac\sigmac}{\zetap-\zetap\zetac} < (>)  \frac{\zetap\Delta + \xip\tau_\mu}{1-\zetap}.
\end{equation*}
For the special case $a/c = b/d$, $F(\chi)\equiv1$ and any $\chi$ is optimal.
\qed

\subsection{Proof of Theorem~\ref{th.3}}\label{ap.th.3}

To impose convergence and performance, we look at the following matrix condition, featuring semidefinite matrix $\X$, the scalar $\rho, \lambda_1, \lambda_2, \gamma_2$, and matrix $\K$ which is implicit in $\B, \B_e$:
\begin{multline}\label{lmicarsten}
(\bullet)^\transp \left[\begin{array}{cc} \X & {\bf 0} \\ {\bf 0} & -\X \end{array} \right]\left[\begin{array}{ccc} {\bf 0} & \B & \B_e \\ \rho \I_n & {\bf 0}_{12} & {\bf 0}_{12} \end{array} \right]  + \lambda_1 (\bullet)^\transp \left[\begin{array}{cc}\omega_1^{2}\I_n & {\bf 0} \\  {\bf 0} & -\I_n \end{array} \right]\left[\begin{array}{cccc} \I_n & {\bf 0} & {\bf 0} & {\bf 0}_{12} \\ {\bf 0} & \I_n & {\bf 0} & {\bf 0}_{12} \end{array} \right]  + \\ + \lambda_2 (\bullet)^\transp \left[\begin{array}{cc}\omega_1^{2}\omega_2^{2}\I_n & {\bf 0} \\ {\bf 0} & -\I_n \end{array} \right]\left[\begin{array}{cccc} \I_n & {\bf 0}  & {\bf 0} & {\bf 0}_{12} \\ {\bf 0} & {\bf 0}  & \I_n & {\bf 0}_{12} \end{array} \right] +\left[\begin{array}{ccc} {\bf 0} & & \\ & {\bf 0}_{22} &  \\ & & - \gamma_2^2\I_n \end{array} \right] \preceq 0,
\end{multline}
where $(\bullet)^\transp$ means that what is post-multiplied is also pre-multiplied transposed and ${\bf 0} = {\bf 0}_{n\times n}$, ${\bf 0}_{ij} = {\bf 0}_{in\times jn}$. Condition~\eqref{lmicarsten} combines the system, the quadratic constraints (i.e., the contractivity) via an $S$-procedure, and the performance criterion. 

We now develop the multiplications, we let $\|\cdot\|^2_{\X}:= (\cdot)^\transp {\X} (\cdot)$, and pre and post multiply with the vector $[(\x_{k}-\hat{\x}^\star_{k+1})^{\transp}, (\bar{\w}_{k+1}-\hat{\x}^\star_{k+1})^{\transp}, (\bar{\uu}_{k+1}-\hat{\x}^\star_{k+1})^{\transp}, {\e}_{k+1}^{\transp} ]^\transp$ and we obtain,
\begin{multline}
-\rho^2\|\x_{k}-\hat{\x}^\star_{k+1}\|_{\X}^2 + \|\x_{k+1}-\hat{\x}^\star_{k+1}\|_{\X}^2  \leq  - \underbrace{\lambda_1 \left[\omega_1^2 \|\x_{k} - \hat{\x}_{k+1}^\star\|^2 - \|\bar{\w}_{k+1} - \hat{\x}_{k+1}^\star\|^2\right]}_{\geq 0} + \\ - \underbrace{\lambda_2 \left[\omega_1^2 \omega_2^2 \|\x_{k} - \hat{\x}_{k+1}^\star\|^2 - \|\bar{\uu}_{k+1} - \hat{\x}_{k+1}^\star\|^2\right]}_{\geq 0} + \gamma^2\| \e_{k+1}\|^2 \leq \gamma^2\| \e_{k+1}\|^2.
\end{multline}

Define the error $E_i := \x_{i}-\hat{\x}^{\star}_{i}$ and the drift $\delta_k = \hat{\x}^{\star}_{k+1}-\hat{\x}^{\star}_{k}$, then 
\begin{equation}
\|E_{k+1}\|^2_{\X} \leq  \rho^2 \|E_{k} - \delta_{k}\|^2_{\X} +  \gamma^2\|\e_{k+1} \|^2. 
\end{equation}
Taking the square root of both sides, since $\geq 0$
\begin{equation}
 \|E_{k+1}\|_{\X} \leq \sqrt{ \rho^2 \|E_{k} - \delta_{k}\|^2_{\X} + \gamma^2\|\e_{k+1} \|^2}  \leq \rho \|E_{k}\|_{\X} + \rho\|\delta_{k}\|_{\X} + \gamma\|\e_{k+1} \|. 
\end{equation}
Let $\|\delta_{k}\| \leq \Delta$, also note that $\E[\|\e_{k+1}\|]\leq 1$ since $\E[\|\q_{k+1}\|]\leq 1/\sqrt{2}, \E[\|\r_{k+1}\|]\leq 1/\sqrt{2}$. Since we impose $\X\succeq \I_n$ without loss of generality (since the problem remains unchanged for any scalar scaling), then $ \|E_{k+1}\|_{\X} \geq  \|E_{k+1}\|$ and $\|E_{k}\|_{\X} \leq \|\X^{1/2}\| \|E_{k}\| = \gamma_1 \|E_{k}\|$. Here the equality sign is due to the fact that we minimize over $\gamma_1$. 

Similarly, $\rho\|\delta_{k}\|_{\X} \leq \rho \gamma_1 \|\delta_{k}\|$. Then, iterating on $k$, and taking the expectations, we obtain,
\begin{eqnarray}\label{ate} 
&&\textbf{E}[\| E_{k}\|] \leq \gamma_1 \rho^{k} \textbf{E}[\|E_{0}\|] + \frac{1}{1-\rho}\left(\gamma_1 \rho \Delta + \gamma_2\right), \\ &&\limsup_{k\to \infty} \textbf{E}[\| E_{k}\|] \leq \frac{1}{1-\rho}\left(\gamma_1 \rho \Delta + \gamma_2\right).
\end{eqnarray}

As such, for any fixed $\rho$, minimizing $\gamma_1 \rho \Delta + \gamma_2$ minimizes the asymptotic tracking error. Furthermore, since $\|x\|_1 \leq \sqrt{n}\|x\|_2$ for $x\in\reals^n$, we know that $\sqrt{\gamma_1^2 \rho^2 \Delta^2 + \gamma_2^2} \geq (\gamma_1 \rho \Delta + \gamma_2)/\sqrt{2}$. So our cost majorizes the asymptotical error $\gamma_1 \rho \Delta + \gamma_2$ and therefore by minimizing our cost, we diminish the latter (notice, we do not minimize the latter, in general, since we have a constrained problem). %

To finish the proof, we need to transform~\eqref{lmicarsten} into~\eqref{condK}. We develop the matrix multiplications, and we observe that the resulting matrix is block diagonal. The first block is $\rho^2 \X \succeq (\lambda_1 \omega_1^2 + \lambda_2 \omega_1^2\omega_2^2) \I_n$. The second block is
\begin{equation}\label{invert}
(\bullet)^\transp \X \left[\B \quad \B_e \right] + \lambda_1 \left[\begin{array}{ccc}-\I_n  & \bf 0 & \bf 0_{12}\\  & \bf 0 & \bf 0_{12} \\ && \bf 0_{22} \end{array} \right] + \lambda_2\left[\begin{array}{ccc}\bf 0  &  \bf 0 & \bf 0_{12}\\  & -\I_n & \bf 0_{12} \\ && \bf 0_{22} \end{array} \right]  + \left[\begin{array}{ccc}\bf 0  & \bf 0 & \bf 0_{12}\\  & \bf 0 & \bf 0_{12} \\ && -\gamma^2 \I_{2n}\end{array} \right] \preceq 0. 
\end{equation}
Expand $\X$ into $\X \X^{-1} \X$ and introduce the variable $\W = \K \X$. Then, taking the Schur's complement, we obtain $\eqref{condK}$, from which the thesis. \qed

\subsection{Proof of Theorem~\ref{th.4}}\label{ap.th.4}

The proof follows the proof of Theorem~\ref{th.3}. We focus here on the different parts. 

Starting from Eq.~\eqref{lmicarsten}, we adapt the first term to: 
\begin{equation}\label{lmicarsten2}
(\bullet)^\transp \left[\begin{array}{cc} \X(\theta_{s+1}) & {\bf 0} \\ {\bf 0} & -\X(\theta_{s}) \end{array} \right]\left[\begin{array}{ccc} {\bf 0} & \B(\theta_{s}) & \B_e(\theta_{s}) \\ \rho \I_n & {\bf 0}_{12} & {\bf 0}_{12} \end{array} \right].
\end{equation}
The bottom diagonal leads to conditions
\begin{eqnarray}
&&\rho^2[\X_0 + \theta_s\X_1] \succeq (\lambda_1 \omega_1^2 + \lambda_2 \omega_1^2  \omega_2^2)\I_n,  \\
&& [\X_0 + \theta_s\X_1] \succeq \I_n, \quad [\X_0 + \theta_s\X_1] \preceq \gamma_1^2 \I_n  
\end{eqnarray}
which needs to be valid for the extreme points $\theta_s=0,1$, since affine in $\theta_s$. However, given the constraint $\X_1 \preceq 0$, the above simplify into~\eqref{eq.affine1} and \eqref{eq.affine2}, respectively. 

Bu using again the constraint $\X_1 \preceq 0$, the upper diagonal can be upper bounded as, 
\begin{multline}\label{lmicarsten3}
(\bullet)^\transp \X(\theta_{s+1})[ \B(\theta_{s}) \quad \B_e(\theta_{s})] = (\bullet)^\transp [\X_0 + \theta_{s+1} \X_1][ \B(\theta_{s}) \quad \B_e(\theta_{s})] \preceq \\ (\bullet)^\transp [\X_0 - \nu \X_1 + \theta_{s} \X_1][ \B(\theta_{s}) \quad \B_e(\theta_{s})] = (\bullet)^\transp {\Y}(\theta_s)[ \B(\theta_{s}) \quad \B_e(\theta_{s})],
\end{multline}
so imposing a $\preceq$ condition on the latter, would imply a condition on the former. In particular, adapting~\eqref{invert}, the condition
\begin{equation}\label{invert2}
(\bullet)^\transp {\Y}(\theta_s)[ \B(\theta_{s}) \quad \B_e(\theta_{s})] + \lambda_1 \left[\begin{array}{ccc}-\I_n  & \bf 0 & \bf 0_{12}\\  & \bf 0 & \bf 0_{12} \\ && \bf 0_{22} \end{array} \right] + \lambda_2\left[\begin{array}{ccc}\bf 0  &  \bf 0 & \bf 0_{12}\\  & -\I_n & \bf 0_{12} \\ && \bf 0_{22} \end{array} \right]  + \left[\begin{array}{ccc}\bf 0  & \bf 0 & \bf 0_{12}\\  & \bf 0 & \bf 0_{12} \\ && -\gamma^2 \I_{2n}\end{array} \right] \preceq 0 
\end{equation}
would imply a similar condition on ${\X}(\theta_{s+1})$, thus the upper diagonal of~\eqref{lmicarsten2}, and for proof of Theorem~\ref{th.3}, convergence of the algorithm as indicated in Theorem~\ref{th.4}.

Condition~\eqref{invert2} leads to condition~\eqref{eq.quadratic}, by Schur complement and dropping the now-redundant subscript $s$.

We remark here the importance of the constraint $\X_1 \preceq 0$, without which we would need two upper bounds in~\eqref{lmicarsten3}, one for $-\nu$ and one for $+\nu$. \rev{This would imply the definition of two matrices $\Y_{-}(\theta) = \X_0 - \nu \X_1 + \theta_s \X_1$ and $\Y_{+}(\theta_s) = \X_0 + \nu \X_1 + \theta_s \X_1$ and the introduction of two auxiliary matrices $\W_{-}(\theta_s) = \K(\theta_s) \Y_{-}(\theta_s)$ and $\W_{+}(\theta_s) = \K(\theta_s) \Y_{+}(\theta_s)$. However, in such a case, since $\W_{-}(\theta_s)$ and $\W_{+}(\theta_s)$ would be considered independent variables, then the determination of $\K(\theta) = \W_{-}(\theta)\Y_{-}(\theta)^{-1} \neq \W_{+}(\theta)\Y_{+}(\theta)^{-1} = \K(\theta)$ is flawed. In fact, without the constraint $\X_1 \preceq 0$, the matrices $\W$'s cannot be introduced and the Schur complement cannot be taken, rendering the overall problem nonlinear in the decision variables. 
 }

\qed

\subsection{Proof of Theorem~\ref{th:MSOS}}

The condition $\theta \in [0,1]$ is equivalent to $\theta (1-\theta) \geq 0$. Then we apply the generalized $S$-procedure as in~\cite{Massioni2020}. 

\qed

\bibliographystyle{IEEEtran}
\bibliography{PaperCollection00}

\end{document}